%% file: uno.tex
\RequirePackage{fix-cm}

\documentclass[smallextended]{svjour3}
\journalname{MPC}

\smartqed

\input{preamble.tex}

\begin{document}


\title{\revision{Implementing a unified solver for nonlinearly constrained optimization}}

\author{Charlie Vanaret \and Sven Leyffer}

\institute{C. Vanaret~\orcidlink{0000-0002-1131-7631} \at
Mathematics and Computer Science Division, 
Argonne National Laboratory, Lemont, IL 60439, USA \\
Department of Mathematical Optimization,
Zuse-Institut Berlin, 14195 Berlin, Germany \\
\email{vanaret@zib.de}
\and
S. Leyffer~\orcidlink{0000-0001-8839-5876} \at
Mathematics and Computer Science Division,
Argonne National Laboratory, Lemont, IL 60439, USA \\
\email{leyffer@anl.gov}
}

\date{\today}

\maketitle

\begin{abstract}
\input{sections/abstract.tex}

\end{abstract}

\input{sections/acknowledgments}
\input{sections/motivation}
\input{sections/notation}
\input{sections/abstract_framework}

\input{sections/state_of_the_art}
\input{sections/uno_description}
\input{sections/combining_the_ingredients}
\input{sections/results}

\input{sections/conclusion}

\bibliographystyle{spmpsci}
\bibliography{includes/robust.bib,includes/NLP.bib,includes/MPEC.bib}

\begin{appendices}
\input{sections/appendices}
\end{appendices}

\end{document}

%% file: preamble.tex
\usepackage{todonotes}
\usepackage{amsfonts,amsmath}
\usepackage[boxruled]{algorithm2e}
\DontPrintSemicolon

\usepackage{subcaption}
\usepackage{xcolor}
\usepackage{bm}
\usepackage[colorlinks=false, urlbordercolor=red]{hyperref}
\usepackage{textpos}
\usepackage{pdflscape}
\usepackage{graphicx}
\usepackage{multirow}
\usepackage{orcidlink}
\usepackage{colortbl}
\usepackage[toc]{appendix}
\usepackage{longtable}
\usepackage{environ}
\usepackage{xfrac}
\usepackage{diagbox}
\usepackage{makecell}

\usepackage{lmodern}
\usepackage{tikz-uml}

\newcommand\equaldef{\stackrel{\text{def}}{=}}

\usepackage{inconsolata}

\graphicspath{{../figures/}}

\newcommand{\R}{\mbox{I}\!\mbox{R}}

\newcommand{\mini}{\mathop{\mbox{min}}}

\newcommand{\st}{\mbox{s.t.}}
\newcommand{\dps}{\displaystyle}

\SetKwComment{Comment}{$\triangleright$\ }{}

\SetCommentSty{commentfont}
\SetKwFunction{FRecurs}{FnRecursive}%
\SetKwProg{Function}{}{}{}\SetKwFunction{FRecurs}{void FnRecursive}%
\SetKwBlock{AcceptanceStep}{Check whether $\hat{x}^{(k+1,l)}$ is acceptable}{end}

\input{includes/ingredients_color_code.tex}

\NewEnviron{globalization_mechanism}{%
	\colorbox{mechanism_color}{
		\parbox{0.97\hsize}{
			\begin{textblock*}{88pt}[1, 1](\hsize,5pt)
			\emph{\scriptsize{globalization mechanism}}
			\end{textblock*}%
			\BODY
		}
	}
}
\NewEnviron{globalization_strategy}{%
	\colorbox{strategy_color}{
		\parbox{0.97\hsize}{
			\begin{textblock*}{75pt}[1, 1](\hsize,5pt)
			\emph{\scriptsize{globalization strategy}}
			\end{textblock*}%
			\BODY
		}
	}
}
\NewEnviron{constraint_relaxation}{%
	\colorbox{constraint_relaxation_color}{
		\parbox{0.97\hsize}{
			\begin{textblock*}{75pt}[1, 1](\hsize,5pt)
			\emph{\scriptsize{constraint relaxation}}
			\end{textblock*}%
			\BODY
		}
	}
}
\NewEnviron{subproblem}{%
	\colorbox{subproblem_color}{
		\parbox{0.97\hsize}{
			\begin{textblock*}{40pt}[1, 1](\hsize,5pt)
			\emph{\scriptsize{subproblem}}
			\end{textblock*}%
			\BODY
		}
	}
}
\newcommand{\highlight}[2][yellow]{\mathchoice%
  {\colorbox{#1}{$\displaystyle#2$}}%
  {\colorbox{#1}{$\textstyle#2$}}%
  {\colorbox{#1}{$\scriptstyle#2$}}%
  {\colorbox{#1}{$\scriptscriptstyle#2$}}}%

\input{includes/notations.tex}

\newcommand{\unoversion}{2.2.0}

\newcommand{\classname}[1]{\texttt{#1}}
\newcommand{\methodname}[1]{\texttt{#1}}

\newcommand{\revision}[1]{\textcolor{black}{#1}}

%% file: includes/ingredients_color_code.tex
\definecolor{constraint_relaxation_color}{HTML}{cab2d6}
\definecolor{inequality_method_color}{HTML}{92b9e1}
\definecolor{mechanism_color}{HTML}{FFCCAC}
\definecolor{strategy_color}{HTML}{FFEB94}
\definecolor{subproblem_color}{HTML}{C1E1DC}

%% file: includes/notations.tex
\newcommand{\constraintmultipliers}{y}

\newcommand{\boundmultipliers}{z}
\newcommand{\Boundmultipliers}{Z}
\newcommand{\barrierparameter}{\mu}
\newcommand{\Lag}{{\cal L}}

\newcommand{\numbervariables}{n}
\newcommand{\numberconstraints}{m}

\newcommand{\trial}[1]{\hat{#1}}

\newcommand{\infeasibilitymeasure}{\eta}
\newcommand{\infeasibilitymodel}{\bm{\infeasibilitymeasure}}

\newcommand{\objectivemeasure}[1]{\omega_{#1}}
\newcommand{\objectivemodel}[1]{\bm{\objectivemeasure{#1}}}

\newcommand{\auxiliarymeasure}{\xi}
\newcommand{\auxiliarymodel}{\bm{\auxiliarymeasure}}

\newcommand{\filterobjectivemeasure}{\phi}
\newcommand{\filterobjectivemodel}{\bm{\filterobjectivemeasure}}

\newcommand{\meritfunction}[1]{\psi_{#1}}
\newcommand{\meritfunctionmodel}[1]{\bm{\meritfunction{#1}}}

\newcommand{\direction}[1]{d_{#1}}

\newcommand{\positiveelastic}{u^+}
\newcommand{\negativeelastic}{u^-}

%% file: sections/abstract.tex
\revision{SQP and interior-point methods (also referred to as Lagrange-Newton methods)} typically share key algorithmic components, such as strategies for computing descent directions and mechanisms that promote global convergence. \revision{Building on this insight}, we introduce a \revision{unifying} framework \revision{with eight building blocks that abstracts the workflows of Lagrange-Newton methods}. We then present Uno, a modular C\texttt{++} solver that implements \revision{our unifying} framework and allows the automatic \revision{combination} of a wide range of strategies
with no programming effort from the user. 
Uno is meant to
(1) organize mathematical optimization strategies into a coherent hierarchy;
(2) offer a wide range of efficient and robust methods that can be compared for a given instance;
(3) enable researchers to experiment with novel optimization strategies; and
(4) reduce the cost of development and maintenance of multiple optimization solvers.
Uno's software design allows user to compose new customized solvers for emerging optimization areas such as robust optimization or optimization problems with complementarity constraints, while building on reliable nonlinear optimization techniques.
We demonstrate that Uno is highly competitive against state-of-the-art solvers filterSQP, IPOPT, SNOPT, MINOS, LANCELOT, LOQO, and CONOPT on a subset of 429 small problems from the \revision{CUTE} collection.
Uno is available as open-source software under the MIT license at \url{https://github.com/cvanaret/Uno} \revision{and via its C, Julia, Python, Fortran, and AMPL interfaces}.

%% file: sections/acknowledgments.tex
\section*{Acknowledgments}
This work was supported by the Applied Mathematics activity within U.S. Department of Energy, Office of Science, Advanced Scientific Computing Research, under Contract DE-AC02-06CH11357.

We would like to thank
David Kiessling (KU Leuven) for fruitful discussions about constrained optimization and his extensive testing of Uno,
Nils-Christian Kempke (Cardinal Operations) for his valuable help with technical questions in C\texttt{++}, and
Gail Pieper (Argonne National Laboratory) for proofreading our manuscript.

\revision{We are highly grateful to the following individuals who contributed their expertise and time to the open-source implementation of Uno: Oscar Dowson, Marcel Jacobse, David Kiessling, Rujia Liu, Stefano Lovato, Alexis Montoison, Manuel Schaich, and Silvio Traversaro. We have greatly benefited from their careful contributions and discussions.}

%% file: sections/motivation.tex
\section{Motivation and contributions}

We consider nonlinearly constrained optimization problems of the form
\begin{equation}
\label{eq:general-NLP}
\begin{array}{ll} \dps
\mini_x  & f(x) \\
\st 	   & l \le
\begin{Bmatrix}
c(x) \\
Ax \\
x
\end{Bmatrix}
\le u,
\end{array}
\end{equation}
where $x \in \R^\numbervariables$, $f : \R^\numbervariables \to \R$, $c: \R^\numbervariables \to \R^{\numberconstraints_c}$, ${A \in \R^{\numberconstraints_A \times \numbervariables}}$, and $l \in (\R \cup \{-\infty\})^{\numberconstraints_c + \numberconstraints_A + \numbervariables}$ and $u \in (\R \cup \{\infty\})^{\numberconstraints_c + \numberconstraints_A + \numbervariables}$.
$f$ and $c$ may be nonconvex, which results in a nonconvex optimization problem.
This formulation allows for unbounded variables and equality constraints and explicitly separates general nonlinear, linear, and bound constraints, enabling solvers to readily exploit this structure. However, for the sake of simplicity of this presentation and without loss of generality, we consider the problem in the following form:
\begin{equation}
\label{eq:NLP}
\tag{NLP}
\begin{array}{ll} \dps
\mini_x  & f(x) \\
\st 	   & c(x) = 0 \\
         & x \ge 0,
\end{array}
\end{equation}
where $x \in \R^\numbervariables$, $f : \R^\numbervariables \to \R$, and $c: \R^\numbervariables \to \R^\numberconstraints$.


Most derivative-based iterative methods for nonlinearly constrained nonconvex optimization (e.g., \cite{GouLey:02,leyffer2010software,nocedalwright2006,fletcher2013practical} share common algorithmic components.
\revision{Building on this insight}, we introduce an abstract framework structured around generic \revision{building blocks} that describes these methods in a unified fashion.
We then present Uno (Unifying \revision{Nonlinear} Optimization)\footnote{Uno was first introduced at the ISMP 2018 conference under the name Argonot~\cite{Vanaret2018Argonot}.}, a modular open-source solver for nonlinearly constrained optimization that unifies \revision{numerous} state-of-the-art methods and organizes existing strategies into a coherent hierarchy.
\revision{Optimization} strategies are implemented as unique software components \revision{(e.g., all line-search methods may be generated based on a single \texttt{LineSearch} class)}; these efficient and robust implementations can be combined at will and compared on a given instance.
This allows users to experiment with new algorithmic ideas by building upon Uno's abstractions and interfaces to modeling languages and subproblem solvers.

Uno \revision{allows the automatic generation} of various strategy combinations on the fly with no programming effort from the user. \revision{While all combinations do not lead to convergent 
methods, some of them result in efficient solvers that may not exist as software implementations.}
We demonstrate that Uno is competitive against state-of-the-art solvers on a subset of 429 \revision{CUTE}~\cite{bongartz1995cute} test problems \revision{translated into AMPL}, while being extensible and lightweight.
We believe that Uno has the potential to serve as an experimental laboratory for practitioners and optimizers
and to accelerate research in \revision{nonlinearly constrained} optimization.
Our ultimate goal is to promote the extension of state-of-the-art nonlinear optimization techniques to new classes of problems such as problems with equilibrium constraints (see, e.g., \cite{OutKocZow:98,LuoPanRal:96,LeyLopNoc:06,Leyf:06,FletLeyf:04}) and nonlinear robust optimization (see, e.g., \cite{Leyfferetal:20}).

This paper is organized as follows. In Section~\ref{sec:notations} we introduce our notation and discuss relevant optimality conditions. In Section~\ref{sec:abstract-framework} we introduce our \revision{unifying} framework with generic \revision{building blocks} for unifying derivative-based iterative methods. In Section~\ref{sec:state-of-the-art} we briefly describe state-of-the-art optimization strategies through the prism of our unifying framework. In Section~\ref{sec:uno} we present the basic algorithmic design of Uno and show how various nonlinear optimization methods fit within the architecture.
In Section~\ref{sec:combining-the-ingredients} we illustrate with \revision{two} concrete strategy combinations how the \revision{building blocks} interact with one another. In Section~\ref{sec:results} we provide preliminary numerical results  and compare Uno against state-of-the-art solvers.

%% file: sections/notation.tex
\section{Notation and stationarity conditions}
\label{sec:notations}

In this section we define our notation and state first-order optimality conditions of \eqref{eq:NLP}.

\subsection{Notation}

We start by defining the scaled Lagrangian or Fritz John function~\cite{burke2014sequential}
of \eqref{eq:NLP} at $(x, \constraintmultipliers, \boundmultipliers, \pi)$:
\begin{equation*}
\Lag_\pi(x, \constraintmultipliers, \boundmultipliers) \equaldef \pi f(x) - \constraintmultipliers^T c(x)  - \boundmultipliers^T x = \pi f(x) - \sum_{j=1}^\numberconstraints \constraintmultipliers_j c_j(x) - \sum_{i=1}^\numbervariables \boundmultipliers_i x_i,
\end{equation*}
where $\constraintmultipliers$ and $\boundmultipliers \ge 0$ are the Lagrange multipliers of the general constraints $c(x) = 0$ and the bound constraints $x \ge 0$, respectively, and $\pi \ge 0$ is an objective multiplier that is introduced to handle infeasibility or lack of constraint qualification (CQ)~\cite{mangasarian1967fritz} in a consistent way.

$\nabla_x \Lag_\pi(x, \constraintmultipliers, \boundmultipliers)$ is the gradient of the scaled Lagrangian with respect to $x$:
\begin{equation*}
\nabla_x \Lag_\pi(x, \constraintmultipliers, \boundmultipliers) \equaldef \pi \nabla f(x) - \sum\limits_{j = 1}^\numberconstraints \constraintmultipliers_j \nabla c_j(x) - \boundmultipliers.
\end{equation*}

$\nabla^2_{xx} \Lag_\pi(x, \constraintmultipliers)$ is the Hessian of the scaled Lagrangian with respect to $x$:
\begin{equation*}
\nabla^2_{xx} \Lag_\pi(x, \constraintmultipliers) = \pi \nabla^2 f(x) - \sum_{j=1}^\numberconstraints \constraintmultipliers_j \nabla^2 c_j(x).
\end{equation*}





\subsection{First-order stationarity conditions}

We are primarily concerned with first-order stationary points. The first-order optimality conditions (aka Fritz John conditions) of problem \eqref{eq:NLP} at a stationary point $x^*$ state that there exist $(\pi^*, \constraintmultipliers^*, \boundmultipliers^*)$ such that
\begin{subequations}
\label{eq:FJ-conditions}
\begin{align}
\begin{split}
\label{eq:stationarity}
\text{(stationarity) } & \quad
\nabla_x \Lag_{\pi^*}(x^*, \constraintmultipliers^*, \boundmultipliers^*) = 0
\end{split} \\
\begin{split}
\label{eq:primal-feasibility}
\text{(primal feasibility) } & \quad
c(x^*) = 0, \quad x^* \ge 0
\end{split}\\
\begin{split}
\label{eq:dual-feasibility}
\text{(dual feasibility) } & \quad
\pi^* \ge 0, \quad \boundmultipliers_i^* \geq 0, \quad (\pi^*, \constraintmultipliers^*, \boundmultipliers^*) \not = (0, 0, 0)
\end{split}\\
\begin{split}
\label{eq:complementarity}
\text{(complementarity) } & \quad
x_i^* \boundmultipliers_i^* = 0, \quad \forall i \in \{ 1, \ldots, \numbervariables\}.
\end{split}
\end{align}
\end{subequations}

If $\pi^* > 0$, the optimality conditions are equivalent to the KKT conditions, which can be recovered by scaling Equation~\eqref{eq:stationarity} by $\sfrac{1}{\pi^*}$.
If $\pi^* = 0$, they characterize Fritz John points, that is, feasible points at which \revision{a} constraint qualification is violated.

%% file: sections/abstract_framework.tex
\section{A \revision{unifying} framework for \revision{nonlinearly constrained} optimization}
\label{sec:abstract-framework}

A local quadratic approximation of \eqref{eq:NLP} \revision{at iteration $k$} about the current \revision{primal-dual} point is given by
\begin{equation}
\begin{array}{lll} \dps
\mini_{x} & \frac{1}{2} (x - x^{(k)})^T H^{(k)} (x - x^{(k)}) + (\nabla f^{(k)})^T (x - x^{(k)}) \\
\st 	& c^{(k)} + (\nabla c^{(k)})^T (x - x^{(k)}) = 0 \\
      & x \ge 0,
\end{array}
\label{eq:subproblem}
\end{equation}
\revision{
where $H(x, \constraintmultipliers)$ is any approximation of $\nabla^2_{xx} \Lag_\pi(x, \constraintmultipliers)$,
$f^{(k)} \equaldef f(x^{(k)})$,
$c^{(k)} \equaldef c(x^{(k)})$,
$\nabla f^{(k)} \equaldef \nabla f(x^{(k)})$,
$\nabla c^{(k)} \equaldef \nabla c(x^{(k)})$, and
$H^{(k)} \equaldef H(x^{(k)}, \constraintmultipliers^{(k)})$.
}
For convex equality-constrained problems, an iteration can be interpreted as taking a Newton step on the first-order optimality conditions of the problem, \revision{hence the name \textit{Lagrange-Newton methods}.}

In this section we introduce a unified view for describing \revision{Lagrange-Newton} methods and argue that they can be assembled by combining \revision{eight generic building blocks or \textit{ingredients}. These ingredients are organized in four layers:}
\begin{itemize}
\item \revision{the \textbf{reformulation} layer:}
   \begin{itemize}
   \item a \colorbox{constraint_relaxation_color}{\textbf{constraint relaxation strategy}} constructs feasible \revision{subproblems} by relaxing the general constraints;
   \item \revision{an \colorbox{inequality_method_color}{\textbf{inequality handling method}} handles the inequality constraints;}
   \end{itemize}
\item \revision{the \colorbox{subproblem_color}{\textbf{subproblem}} layer:}
   \begin{itemize}
   \item \revision{a \textbf{Hessian model} determines the approximation of the Lagrangian Hessian;}
   \item \revision{an \textbf{inertia correction strategy} corrects the inertia of the Lagrangian Hessian or of the KKT matrix;}
   \end{itemize}
\item \revision{the \textbf{subproblem solver} layer:}
   \begin{itemize}
   \item \revision{a \textbf{subproblem solver} approximately solves the subproblem by exploiting its characteristics, such as convexity, presence of constraints, presence of bounds, and presence of inequality constraints;}
   \end{itemize}
\item \revision{the \textbf{globalization} layer:}
   \begin{itemize}
   \item a \colorbox{strategy_color}{\textbf{globalization strategy}} determines whether a trial iterate makes sufficient progress toward a solution and accepts or rejects it; and,
   \item a \colorbox{mechanism_color}{\textbf{globalization mechanism}} \revision{controls the length of the direction and} defines the recourse action taken when a trial iterate is rejected.
   \end{itemize}
\end{itemize}
This coloring will be used throughout to illustrate how \revision{the layers and} ingredients interact with one another within an optimization algorithm.
The role of each ingredient is shown in the following abstract algorithm (Algorithm~\ref{alg:iterative-methods}).

\begin{algorithm}
\SetAlgoVlined
\caption{\revision{Abstract algorithm.}}
\small
\KwData{initial point $(x^{(0)}, \constraintmultipliers^{(0)}, \boundmultipliers^{(0)})$}
Set $k \gets 0$ \;
\While{termination criteria at $(x^{(k)}, \constraintmultipliers^{(k)}, \boundmultipliers^{(k)})$ not met}{
    \begin{globalization_mechanism}
    \Repeat{$(\trial{x}^{(k+1)}, \trial{\constraintmultipliers}^{(k+1)}, \trial{\boundmultipliers}^{(k+1)})$ is \colorbox{strategy_color}{acceptable}}{
        Solve a (sequence of) \colorbox{constraint_relaxation_color}{feasible} \colorbox{subproblem_color}{subproblem}(s) that approximate(s) \eqref{eq:NLP} at $(x^{(k)}, \constraintmultipliers^{(k)}, \boundmultipliers^{(k)})$ \;
        Assemble trial iterate $(\trial{x}^{(k+1)}, \trial{\constraintmultipliers}^{(k+1)}, \trial{\boundmultipliers}^{(k+1)})$ \;
    }
    \end{globalization_mechanism}
    Update $(x^{(k+1)}, \constraintmultipliers^{(k+1)}, \boundmultipliers^{(k+1)}) \gets (\trial{x}^{(k)}, \trial{\constraintmultipliers}^{(k)}, \trial{\boundmultipliers}^{(k)})$ \;
    $k \gets k+1$ \;
} 
\KwResult{$(x^{(k)}, \constraintmultipliers^{(k)}, \boundmultipliers^{(k)})$}
\label{alg:iterative-methods}
\end{algorithm}

The inner loop (\textbf{repeat}) generates and solves a feasible subproblem (possibly a sequence of feasible subproblems) until a trial iterate is accepted by the globalization strategy, and the outer loop (\textbf{while}) generates a sequence of acceptable iterates until termination.

\begin{figure}[htbp]
\centering
\includegraphics[width=\columnwidth]{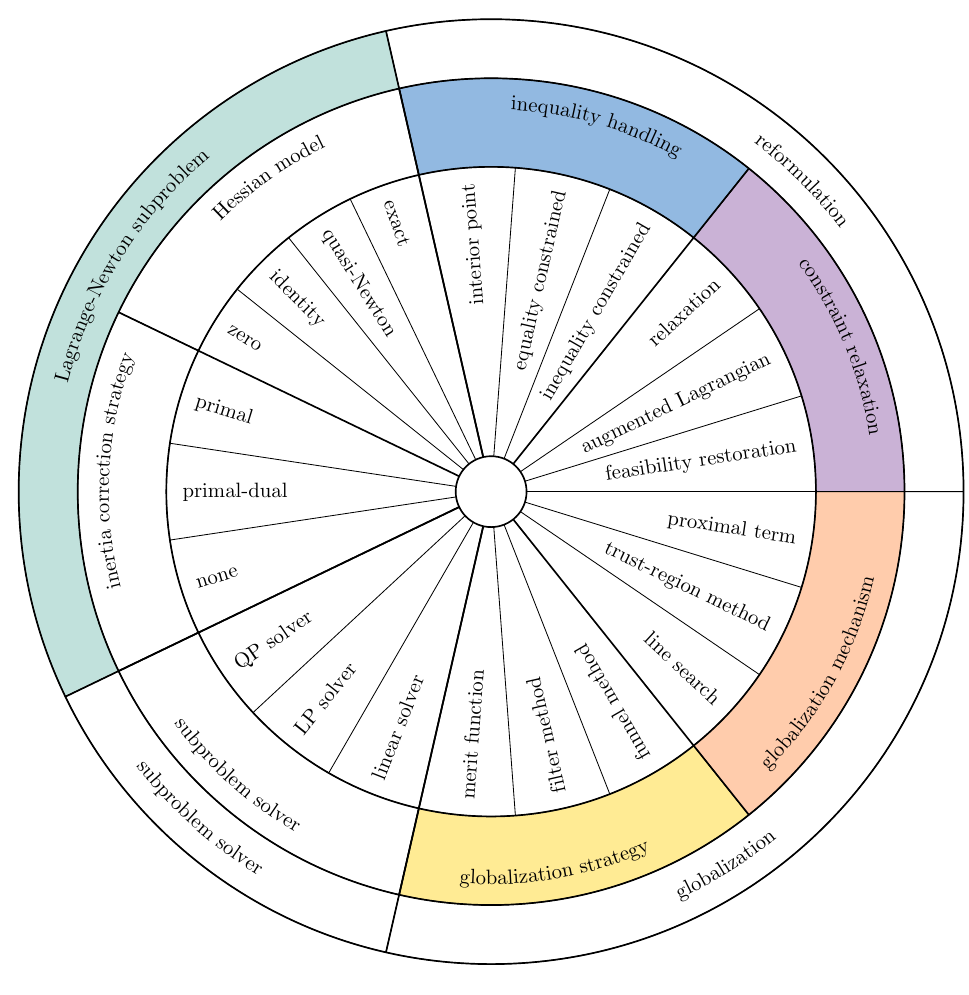}
\caption{\revision{Unifying framework: the ``wheel of strategies''.}}
\label{fig:wheel-of-strategies}
\end{figure}

\revision{
The ``wheel of strategies'' (Figure~\ref{fig:wheel-of-strategies}) organizes a wide range of strategies into a coherent hierarchy.
The outer layer represents the optimization layers, the middle layer represents the ingredients, and the inner layer represents the strategies.
Strategies that perform similar tasks within an optimization method are listed under the same ingredient (e.g., a line search and a trust-region method are both globalization mechanisms).
}

%% file: sections/state_of_the_art.tex
\section{Unified view of state-of-the-art techniques}
\label{sec:state-of-the-art}

In this section we describe a set of strategies that fall into each of the \revision{following} ingredients of our \revision{unifying} framework: \revision{globalization strategies, constraint relaxation strategies, inequality handling methods, Hessian models, inertia correction strategies, and globalization mechanisms}. Our goal is to illustrate the wide variety of methods within a common notation to motivate the design of our modular solver.

\subsection{Globalization strategies}
\label{sec:globalization-strategies}

Constrained optimization methods must achieve two competing goals: minimizing the objective function and minimizing the constraint violation.
Globalization strategies determine whether a trial iterate $\trial{x}^{(k+1)} \equaldef x^{(k)} + \alpha \direction{x}^{(k)}$ (\revision{for a given step length} $\alpha \in (0, 1]$) makes acceptable progress with respect to these goals.
We consider strategies that ensure \textit{global convergence}, that is, convergence to a local minimum, or (weaker) stationary point, from any starting point.
In addition, ideally, the minimization of the measure of infeasibility takes precedence. However, in problems where no feasible point exists, a (global) minimum of the constraint violation is a certificate that the problem is infeasible.

Three (possibly primal-dual) \textit{progress measures} are monitored throughout the optimization process:
\begin{enumerate}
\item an \textbf{infeasibility measure} $\infeasibilitymeasure$ (typically $\|c(x)\|$ for some norm);
\item an \textbf{objective measure} $\objectivemeasure{\pi}$ parameterized by the objective multiplier $\pi \ge 0$ (typically $\pi f(x)$); and
\item an \textbf{auxiliary measure} $\auxiliarymeasure$ (terms such as barrier and proximal terms).
\end{enumerate}

The local models of $\infeasibilitymeasure(x)$, $\objectivemeasure{\pi}(x)$, and $\auxiliarymeasure(x)$ at iteration $k$ about the current iterate are denoted by $\infeasibilitymodel^{(k)}(\direction{x})$, $\objectivemodel{\pi}^{(k)}(\direction{x})$, and $\auxiliarymodel^{(k)}(\direction{x})$ for a given primal direction $\direction{x}$. We define the respective \textit{predicted reductions}
\begin{equation*}
\begin{aligned}
\Delta \infeasibilitymodel^{(k)}(\direction{x}) 	& \equaldef \infeasibilitymodel^{(k)}(0) - \infeasibilitymodel^{(k)}(\direction{x}) \\
\Delta \objectivemodel{\pi}^{(k)}(\direction{x}) 	& \equaldef \objectivemodel{\pi}^{(k)}(0) - \objectivemodel{\pi}^{(k)}(\direction{x}) \\
\Delta \auxiliarymodel^{(k)}(\direction{x}) 			& \equaldef \auxiliarymodel^{(k)}(0) - \auxiliarymodel^{(k)}(\direction{x})
\end{aligned}
\end{equation*}



In order to ensure convergence, \revision{the progress measures and their local models} must be intimately linked to the \revision{definitions of the problem and the subproblem (see Section~\ref{sec:uno-progress-measures})}.

\begin{figure}[htbp!]
\centering
\includegraphics[width=0.9\columnwidth]{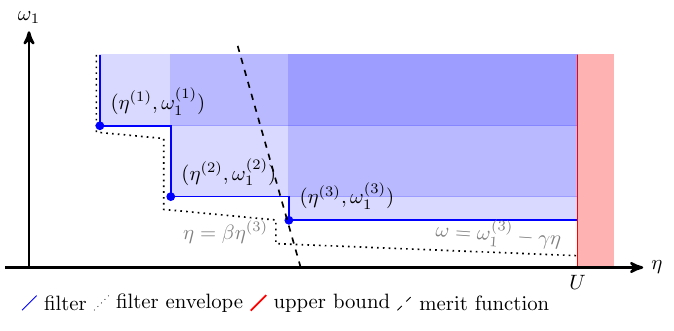}
\caption{Example of a filter and a merit function.}
\label{fig:strategies}
\end{figure}

The two main classes of globalization strategies are merit functions and filter \revision{methods} (Figure \ref{fig:strategies}).
They typically enforce a sufficient decrease condition that forces some scalar combination of $\infeasibilitymeasure$, $\objectivemeasure{\pi}$, and $\auxiliarymeasure$ to decrease by at least a fraction of the decrease predicted by the local model.

\subsubsection{Merit functions}

A merit (or penalty) function combines the three goals $\infeasibilitymeasure$, $\objectivemeasure{\rho}$, and $\auxiliarymeasure$ into a single scalar value:
\begin{equation*}
\meritfunction{\rho}(x) \equaldef \objectivemeasure{\rho}(x) + \infeasibilitymeasure(x) + \auxiliarymeasure(x),
\end{equation*}
where $\rho \ge 0$ is an inverse penalty parameter.
Its predicted reduction is given by
\begin{equation*}
\Delta \meritfunctionmodel{\rho}^{(k)}(\direction{x}) \equaldef \Delta \objectivemodel{\rho}^{(k)}(\direction{x}) + \Delta \auxiliarymodel^{(k)}(\direction{x}) + \Delta \infeasibilitymodel^{(k)}(\direction{x}).
\end{equation*}
The trial iterate $\trial{x}^{(k+1)}$ is accepted if the actual reduction $\meritfunction{\rho}(x^{(k)}) - \meritfunction{\rho}(\trial{x}^{(k+1)})$ in the merit function is larger than a fraction $\sigma \in (0, 1)$ of its predicted reduction:
\begin{equation}
\label{eq:merit-decrease}
\meritfunction{\rho}(x^{(k)}) - \meritfunction{\rho}(\trial{x}^{(k+1)}) \ge \sigma \Delta \meritfunctionmodel{\rho}^{(k)}(\alpha \direction{x}^{(k)}).
\end{equation}
The sufficient decrease condition \eqref{eq:merit-decrease} (also known as Armijo condition) indicates how the subproblem solve is connected to the merit function to ensure global convergence.

\subsubsection{Filter methods}

Filter methods are motivated by the desire to decouple reduction in the objective function from progress toward feasibility. We can interpret a filter method as a mechanism to force iterates closer to the feasible region, so that unconstrained sufficient reduction conditions on the objective can be used to force convergence (see, e.g., \cite{FleLeyToi:02,Fletcher2002Global,WaeBie:05a,WaeBie:05b}).
The decrease function is given by
\begin{equation*}
\filterobjectivemeasure(x) \equaldef \objectivemeasure{1}(x) + \auxiliarymeasure(x),
\end{equation*}
and its predicted reduction is given by
\begin{equation*}
\Delta \filterobjectivemodel^{(k)}(\direction{x}) \equaldef \Delta \objectivemodel{1}^{(k)}(\direction{x}) + \Delta \auxiliarymodel^{(k)}(\direction{x}).
\end{equation*}

Filter methods measure progress toward a solution by comparing the trial infeasibility measure $\infeasibilitymeasure$ and objective measure $\filterobjectivemeasure$ to a filter $\mathcal{F}$, a list of pairs $(\infeasibilitymeasure^{(l)}, \filterobjectivemeasure^{(l)})$ (typically from previous iterates) such that no pair dominates another pair; that is,  there exists no index $l'$ such that $\infeasibilitymeasure^{(l')} < \infeasibilitymeasure^{(l)}$ and $\filterobjectivemeasure^{(l')} < \filterobjectivemeasure^{(l)}$
for all $(\infeasibilitymeasure^{(l)}, \filterobjectivemeasure^{(l)}) \in \mathcal{F}$. Formally, the trial iterate $\trial{x}^{(k+1)}$ is acceptable to the filter if and only if the following conditions hold:
\begin{equation*}
\filterobjectivemeasure(\trial{x}^{(k+1)}) \leq \filterobjectivemeasure^{(l)} - \gamma \infeasibilitymeasure(\trial{x}^{(k+1)}) \quad
\mbox{ or } \quad
\infeasibilitymeasure(\trial{x}^{(k+1)}) < \beta \infeasibilitymeasure^{(l)}, \;  \quad \forall (\infeasibilitymeasure^{(l)}, \filterobjectivemeasure^{(l)}) \in {\cal F}^{(k)}, 
\end{equation*}
where $\gamma > 0$ and $0 < \beta < 1$ are constants that ensure that iterates cannot accumulate at infeasible limit points. These conditions are represented as the filter envelope in Figure~\ref{fig:strategies}.

The filter provides convergence only to a feasible limit because any infinite sequence of iterates must converge to a point where $\infeasibilitymeasure(x) = 0$, provided that $\filterobjectivemeasure(x)$ is bounded below. To ensure convergence to a local minimum, filter methods use a standard sufficient reduction (Armijo) condition from unconstrained optimization:
\begin{equation}
\label{eq:filter-reduction}
\filterobjectivemeasure(x^{(k)}) - \filterobjectivemeasure(\trial{x}^{(k+1)}) \geq \sigma \Delta \filterobjectivemodel^{(k)}(\alpha \direction{x}^{(k)}),
\end{equation}
where $\sigma \in (0, 1)$.
It makes sense to enforce this condition only if the model predicts a decrease in the objective function.
Thus, filter methods use the switching condition
\begin{equation*}
\Delta \filterobjectivemodel^{(k)}(\alpha \direction{x}^{(k)}) \geq \delta \infeasibilitymeasure(x^{(k)})^2 ,
\end{equation*}
where $\delta > 0$, to decide when \eqref{eq:filter-reduction} should be enforced.
If the trial point is accepted, it is added to $\mathcal{F}^{(k)}$ if $\infeasibilitymeasure(x^{(k)}) > 0$ or if the switching condition fails (which automatically satisfies $\infeasibilitymeasure(x^{(k)}) > 0$).

\subsection{Constraint relaxation strategies}


\revision{Subproblem constraint relaxation strategies, referred to as constraint relaxation strategies for the sake of conciseness, guarantee that the subproblems are well defined and feasible. }
In general, we cannot assume that the nonlinear problem \eqref{eq:NLP} is \revision{locally or globally} feasible, \revision{nor that a given subproblem is feasible (the linearized constraints may be inconsistent or the trust-region radius may be too small)}.
Moreover, nonlinear solvers may converge to points that violate standard constraint qualifications, and we must take these situations into account when defining optimality conditions.
\revision{To address these unfavorable scenarios, constraint relaxation strategies relax the constraints systematically or when the standard subproblem is infeasible, thus giving more emphasis to decreasing the constraint violation.}
We review two constraint relaxation strategies: $\ell_1$ relaxation and feasibility restoration.

\subsubsection{$\ell_1$ relaxation}
\label{sec:l1-relaxation}

The $\ell_1$ relaxation strategy replaces a constrained optimization problem with a nonsmooth bound-constrained problem in which a penalty term is added to the objective:
\begin{equation}
\begin{array}{ll} \dps
\mini_x	& \rho f(x) + \|c(x)\|_1 \\
\st     & x \ge 0,
\end{array}
\label{eq:l1-relaxation}
\end{equation}
where $\rho \ge 0$ is an inverse penalty parameter. An appropriate value of $\rho$ \revision{may be determined by steering rules}~\cite{Byrd2008Steering,Byrd2010Infeasibility,burke2014sequential}.

\revision{Note that the penalty parameter is often attached to \revision{$\|c(x)\|_1$} in the literature. We adopt the inverse notation as in \cite{Byrd2010Infeasibility} for several reasons: (i) it is numerically easier to drive $\rho$ to 0 than to drive the penalty parameter of $\|c(x)\|_1$ to $+\infty$; (ii) in second-order methods, the inverse penalty parameter enters the Hessian as objective multiplier; and (iii) as $\rho \rightarrow 0$, only $\objectivemeasure{\rho}$ vanishes, and the implicit constraints in $\auxiliarymeasure$ (barrier or proximal terms) are still enforced.}

The nonsmooth $\ell_1$ relaxed problem can be reformulated as a smooth constrained problem with elastic variables.
Other norms can be used; however, the $\ell_1$ norm has emerged as the preferred option. In particular, the $\ell_1$ relaxation is exact: one can show under mild conditions that for $\rho > 0$ sufficiently small,
a second-order sufficient point of \eqref{eq:l1-relaxation} is also a second-order sufficient point of \eqref{eq:NLP}
(see, e.g., Theorem 14.3.1 in~\cite{fletcher2013practical}).

\subsubsection{Feasibility restoration}
\label{sec:feasibility-restoration}

An infeasible subproblem results from inconsistent linearized or bound constraints, \revision{or when the trust-region radius is too small}; it is an indication that \eqref{eq:NLP} may be infeasible. In this case, the method reverts to the \textit{feasibility restoration phase}: the original objective is temporarily discarded, and the following feasibility problem is solved instead:
\begin{equation}
\begin{array}{ll} \dps
\mini_x & \| c(x) \| \\
\st     & x \ge 0,
\end{array}
\label{eq:feasibility-problem}
\end{equation}
for some norm $\| \cdot \|$ in $\R^\numberconstraints$.
The aim of solving the feasibility problem is to compute a point as close as possible to the feasible region. Feasibility restoration improves feasibility until a minimum of the constraint violation is obtained or the subproblem becomes feasible again, in which case solving the original problem (the \textit{optimality phase}) is resumed.
Any (local) solution $x^* \geq 0$ of \eqref{eq:feasibility-problem} with $\| c(x^*) \| > 0$ is a certificate that \eqref{eq:NLP} is (locally) infeasible.

\subsection{\revision{Inequality handling methods}}

\revision{Inequality handling methods deal with the combinatorics induced by inequality constraints within the context of Lagrange-Newton approximations.}
We briefly review three classes of \revision{methods}: inequality-constrained methods, equality-constrained methods, and interior-point methods.
Note that a classification into inequality-constrained and equality-constrained SQP methods can be found in~\cite{nocedalwright2006}.

\subsubsection{Inequality-constrained \revision{methods}}

In inequality-constrained methods, the handling of inequality constraints is deferred to the subproblem solver. Traditionally, sequential quadratic programming (SQP) methods~\cite{WilRB:63,HanSP:77,PowMJD:78a,PowYua:86} generate a sequence of \revision{inequality-constrained} quadratic problems (QPs) that are solved by means of an active-set QP solver: it maintains an estimate of the active set and solves a sequence of equality-constrained subproblems. The estimate of the active set is updated at each iteration using dual information, until the algorithm terminates with primal-dual feasibility. SQP methods converge quadratically under reasonable assumptions near a local minimizer, once the active set settles down.
Sequential linear programming~\cite{ChiFle:03} is a particular case of SQP in which the subproblems are linear problems (LPs); that is, no second-order information is exploited ($H^{(k)} = 0$).

\subsubsection{Equality-constrained \revision{methods}}

Equality-constrained methods operate in two phases: the first phase solves a ``low-fidelity'' subproblem (such as an LP or a convex QP with a quasi-Newton Hessian), which provides an estimate of the active set $\mathcal{A}$. The second phase solves a ``high-fidelity'' (such as second-order) equality-constrained problem with exact Hessian in which the inequality constraints of the active set are fixed to their active bounds and the inactive inequalities are dropped. This is illustrated in the following equality-constrained problem (using the notation of Problem~\eqref{eq:general-NLP}):
\begin{equation*}
\begin{array}{ll} \dps
\mini_x  & f(x) \\
\st 	   & 
\begin{Bmatrix}
c(x) \\
Ax \\
x
\end{Bmatrix}_{\mathcal{A}}
= b_{\mathcal{A}},
\end{array}
\end{equation*}
where each component of $b_{\mathcal{A}}$ is the active (lower or upper) bound of the corresponding constraint.
A typical example of equality-constrained methods is SLPEQP~\cite{ChiFle:03} (aka SLQP): it estimates the active set by solving an LP, then solves an equality-constrained QP (EQP).

\subsubsection{Interior-point \revision{methods}}
\label{sec:ipm}

\textit{Primal-dual interior-point methods} relax the complementarity equations \eqref{eq:complementarity}
by a \revision{barrier parameter} $\barrierparameter > 0$:
\begin{equation*}
x_i \boundmultipliers_i = \barrierparameter, \quad \forall i \in \{ 1, \ldots, \numbervariables\},
\end{equation*}
which implies $x > 0$ and $\boundmultipliers > 0$.
Positivity of $x$ and $\boundmultipliers$ is then enforced at each iteration, \revision{e.g,} by a fraction-to-the-boundary rule~\cite{Waechter2006Implementation}.
Similarly to homotopy methods, $\barrierparameter$ is asymptotically driven to 0 until the termination criteria are met.

In contrast, \revision{\textit{primal interior-point methods} (or barrier methods) treat $\boundmultipliers$ as dependent variables. \revision{This corresponds to solving the following nonlinear equality-constrained problem}:
\begin{equation*}
\begin{array}{ll} \dps
\mini_x  & f(x) - \mu \dps \sum_{i=1}^n \log(x_i) \\
\st 	   & c(x) = 0.
\end{array}
\end{equation*}
}

\subsection{\revision{Hessian models}}
\label{sec:hessian-models}

\revision{
The Hessian model is a high- or low-quality representation of the Lagrangian Hessian of the original model. It may be:
\begin{itemize}
\item the exact Hessian $H^{(k)} = \nabla_{xx}^2 \mathcal{L}_\rho(x^{(k)}, \constraintmultipliers^{(k)})$ if it is provided by the modeler or the modeling framework;
\item a quasi-Newton approximation such as BFGS or SR1, or a limited-memory version thereof;
\item the identity matrix;
\item the zero matrix.
\end{itemize}
As the model gets reformulated by the constraint relaxation strategy and the inequality handling method, structured (typically diagonal) curvature may be introduced in the form of a proximal term or a barrier term.
}

\subsection{\revision{Inertia correction strategies}}

\revision{
Since line-search methods require that $\direction{x}^{(k)}$ be a descent direction, the Lagrangian Hessian or the augmented system is usually regularized such that the former becomes positive definite, and the latter has $n$ positive, $m$ negative, and $0$ zero eigenvalues, that is the inertia $(n, m, 0)$. An approximation of the inertia is provided by factorization codes such as MA57~\cite{Duff2004MA57}.
}

\revision{
Typically, the inertia of the Hessian $H^{(k)}$ is corrected by adding a multiple of the identity matrix:
\begin{equation*}
H^{(k)} + \delta I
\end{equation*}
for some ${\delta > \min\{ \text{eig}( H^{(k)}) \}}$.
The inertia of an augmented matrix
\begin{equation}
M \equaldef
\begin{pmatrix}
A & B \\
B^T & \quad \\
\end{pmatrix}
\end{equation}
is controlled by adjusting the regularization parameters $\delta_w \ge 0$ and $\delta_c \ge 0$ in the matrix:
\begin{equation}
\begin{pmatrix}
A + \delta_w I & B \\
B^T & \quad -\delta_c I \\
\end{pmatrix}.
\end{equation}
}


\subsection{Globalization mechanisms}


When the methods are started far from a solution, directions may be unbounded or result in trial iterates that increase both the objective and the constraint violation.
Globalization mechanisms provide a recourse action if a local approximation is deemed too poor to make progress toward a solution:
line-search methods restrict the length of the step along a given direction,
while
trust-region methods restrict the length of the direction a priori.

\subsubsection{Line-search methods}

Line-search methods compute a trial step $\alpha \direction{x}^{(k)}$ by determining a step length $\alpha \in (0, 1]$ along the direction $\direction{x}^{(k)}$.
Inexact line-search methods determine an approximate step length \revision{such that the corresponding trial iterate is accepted} by the globalization strategy (sufficient decrease for a merit function or acceptance by a filter).
For instance, backtracking line-search methods generate a sequence of trial iterates $x^{(k)} + \alpha^{(k,l)} \direction{x}^{(k)}$ for $l \in \mathbb{N}$ until acceptance, where ${\alpha^{(k,l)} = c^l \alpha_{max}^{(k)}}$, $\alpha_{max}^{(k)} \in (0, 1]$ (usually $1$ in SQP methods and smaller than $1$ in IPM) and $c \in (0, 1)$.


Line-search methods require that the approximation $H^{(k)}$ of the Lagrangian Hessian be positive definite on the nullspace of the Jacobian of the active constraints in order to guarantee that $\direction{x}^{(k)}$ is a descent direction for the
merit function.

\subsubsection{Trust-region methods}

Trust-region methods limit the length of the direction $\direction{x}$ a priori by imposing the trust-region constraint $\|\direction{x}\| \leq \Delta^{(l)}$ for some norm $\| \cdot \|$ in $\R^\numbervariables$, where $\Delta^{(l)} > 0$ is the trust-region radius. The step $\direction{x}^{(k, l)}$ is then computed by (approximately) solving the trust-region subproblem.
If the trial iterate ${\trial{x}^{(k+1, l)} \equaldef x^{(k)} + \direction{x}^{(k, l)}}$ is accepted by the globalization strategy, $\Delta^{(l)}$ is increased if the trust region is active at $\direction{x}^{(k, l)}$. Otherwise, $\Delta^{(l)}$ is decreased to a value smaller than $\min(\Delta^{(l)}, \lVert \direction{x}^{(k,l)} \rVert)$.
Contrary to line-search methods, a positive definite Lagrangian Hessian is not required because directions of negative curvature are bounded by the trust region.

\subsection{Summary}

Table~\ref{tab:state-of-the-art-solvers} presents a unified view of state-of-the-art solvers ALGLIB MinNLC~\cite{bochkanov2011alglib}, CONOPT~\cite{drud1994conopt}, FICO XSLP~\cite{FICO2023}, filterSQP~\cite{FleLey:98}, IPOPT~\cite{Waechter2006Implementation}, KNITRO~\cite{ByrdNoceWalt:04}, LANCELOT~\cite{ConGouToi:92}, LOQO~\cite{vanderbei1999loqo}, MINOS~\cite{MurSaun:93}, NAG e04uc/e04wdc~\cite{NAG2017}, NLPQL~\cite{schittkowski2001nlpqlp}, SLSQP~\cite{kraft1988software}, SNOPT~\cite{GilMurSau:02}, SQuID~\cite{burke2014sequential}, and WORHP~\cite{buskens2013esa}.
Each solver is characterized in terms of the ingredients within the proposed \revision{unifying} framework.
Note that FICO XSLP is the only solver that does not implement a proper globalization strategy.


%% file: sections/uno_description.tex
\section{Uno: a \revision{unified solver for nonlinearly constrained optimization}}
\label{sec:uno}

We have implemented our unifying framework for nonlinearly constrained optimization within Uno, a modular solver written in C\texttt{++}17. A generic and flexible code, it supports a broad range of strategies that can be combined automatically and on the fly with no programming effort from the user.
Additional software components are implemented \revision{orthogonally to the \revision{eight} ingredients}, such as termination criteria.
The code is packaged in a lightweight library (around \revision{10,000} lines of code for the current version Uno \revision{\unoversion}, \revision{excluding interfaces}) 
available as open-source software under the MIT license at \url{https://github.com/cvanaret/Uno}~.
Uno \revision{\unoversion} is \revision{available via its C, Julia (registered package \texttt{UnoSolver.jl}), Python (the packaging of \texttt{unopy} is under development), Fortran, and AMPL~\cite{FouGayKer:03,Gay1997Hooking}} interfaces.

In the following, we briefly present the generic architecture of Uno \revision{\unoversion}, describe how ingredients can be automatically combined, list its features, and discuss the current limitations.


\subsection{Generic architecture}

The modularity of Uno stems from its generic architecture: each ingredient is implemented \revision{once and} independently of the others, which
improves readability, and makes building blocks less prone to error and easier to maintain than monolithic codes. This results in a modern, flexible, and maintainable framework.


Figure~\ref{fig:uno-uml} represents Uno's object-oriented architecture as a Unified Modeling Language (UML)\footnote{\url{https://www.visual-paradigm.com/guide/uml-unified-modeling-language/uml-class-diagram-tutorial/}}) diagram based on inheritance (``is a'') and composition (``has a'').
The ingredients are modeled as abstract classes: they define interfaces, that is, generic actions and behaviors that must be implemented by concrete strategies modeled as subclasses.
For example:
\begin{itemize}
\item the classes \classname{BacktrackingLineSearch} and \classname{TrustRegionStrategy} inherit \revision{from} the abstract class \classname{GlobalizationMechanism} (inheritance is represented by dashed arrows) and thus must implement the purely virtual member function \methodname{compute\_acceptable\_iterate()};
\item the \classname{GlobalizationMechanism} class possesses a member of type \\\classname{ConstraintRelaxationStrategy} (composition is represented by solid diamond lines).
\end{itemize}


\subsection{Automatic strategy combinations}
\label{sec:combinations}

Uno \revision{\unoversion} implements state-of-the-art strategies that can be combined automatically thanks to the modular software architecture. The number of possible strategy combinations is the size of the Cartesian product of the \revision{eight} ingredients.
Note that all combinations do not necessarily result in sensible algorithms, or even convergent approaches.



At the moment, Uno prohibits the combination \revision{of \colorbox{subproblem_color}{interior-point methods} and \colorbox{mechanism_color}{trust-region strategies}}. A possible \revision{strategy is KNITRO's step decomposition: the direction is decomposed into a normal step that minimizes the constraint violation within the trust region, and a tangential step that minimizes the objective for a given constraint violation target within the trust region}.
This limitation will be resolved in later Uno versions.

Some strategy combinations that correspond to state-of-the-art solvers
are available as ``presets'': \revision{they automatically connect} the \revision{eight} ingredients, and \revision{set values for} the hyperparameters that can be found in the solvers' documentations (they will be listed in the Uno user manual). The following presets are available in Uno \revision{\unoversion}:
\begin{itemize}
\item \texttt{filtersqp}: A trust-region restoration filter SQP method \`a la filterSQP~\cite{Fletcher2002Global}. Second-order correction steps were not implemented.
\item \texttt{ipopt}: A line-search restoration filter interior-point method \`a la IPOPT~\cite{Waechter2006Implementation}. Second-order correction steps, \revision{scaling, least-square multipliers}, iterative refinement, iterative bound relaxations, non-monotone techniques, and soft feasibility restoration were not implemented.
\end{itemize}

\subsection{\revision{Features of Uno \unoversion}}

Uno \revision{\unoversion} contains a number of features that are described here.

\subsubsection{Interfaces to subproblem solvers}

Interfaces to the following subproblem solvers are available: 
\begin{itemize}
\item BQPD~\cite{FleR:95,fletcher2000stable}: a null-space active-set solver for nonconvex QPs. \revision{BQPD accepts Hessian-vector products instead of an explicit matrix. A Hessian linear operator can be provided to the Uno model};
\item MA57~\cite{Duff2004MA57}, \revision{MA27~\cite{duff1991factorization}, and MUMPS~\cite{amestoy2000mumps}}: direct solvers for sparse symmetric indefinite linear systems;
\item \revision{HiGHS~\cite{huangfu2018parallelizing}: a parallel simplex implementation for linear programming. Uno currently does not support the HiGHS QP solvers}.
\end{itemize}

\revision{
We plan to write an interface to GALAHAD's TRS package~\cite{fowkes2023galahad} that solves the trust-region problem to global optimality:
\begin{equation*}
\begin{array}{ll} \dps
\mini_d     & \frac{1}{2} d^T H d + c^T d + f \\
\st 	    & Ad = 0 \\
            & \|d\|_2 \le \Delta.
\end{array}
\end{equation*}
}

\subsubsection{\revision{Definition of the subproblem}}
\label{sec:uno-subproblem}

\revision{
The constraint relaxation strategy and the inequality handling method successively reformulate the original problem with respect to the general constraints $c(x) = 0$ and the bound constraints $x \ge 0$. A local model of the resulting reformulated problem is then built.
}

\revision{
The Lagrange-Newton subproblem is composed of the following elements: the reformulated problem, the current primal-dual iterate, the Hessian model, the inertia correction strategy, and a possible trust-region radius. These elements interact with one another to automatically define the progress measures and their local models (see Section~\ref{sec:uno-progress-measures}), the Lagrangian Hessian or augmented matrix with the correct inertia (see Section~\ref{sec:uno-inertia-correction}), and the type of subproblem solvers that can be used (see Section~\ref{sec:instantation-ingredients}).
}

\begin{figure}[h!]
\centering
\includegraphics[width=\columnwidth]{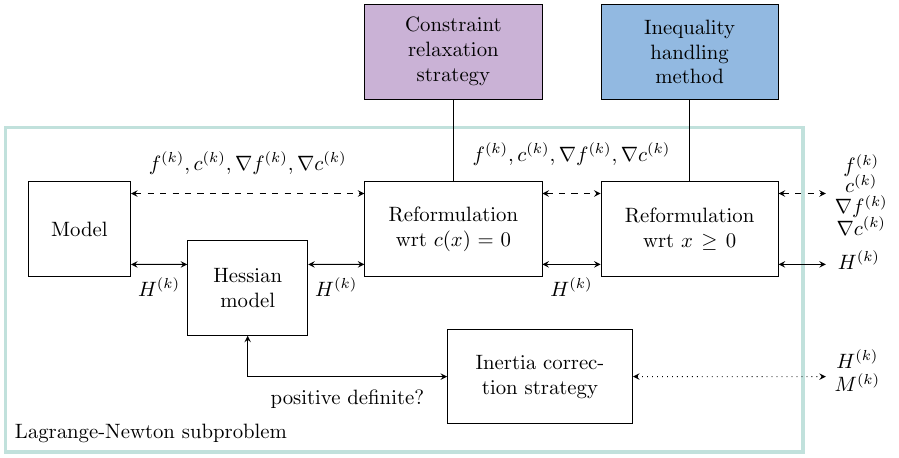}
\caption{\revision{Definition of a \texttt{Subproblem}.}}
\label{fig:uno-subproblem}
\end{figure}

\revision{
The subproblem is represented schematically in Figure~\ref{fig:uno-subproblem}. Dashed arrows represent zeroth- and first-order calls. Full arrows represent calls to the Lagrangian Hessian. Dotted arrows represent calls to the inertia correction strategy for a given Hessian $H^{(k)}$ or augmented system $M^{(k)}$.
}

\subsubsection{\revision{Definition of progress measures}}
\label{sec:uno-progress-measures}

\revision{
The progress measures (see Section~\ref{sec:globalization-strategies}) are defined automatically by the successive reformulations of the optimization problem (\texttt{l1RelaxedProblem} and \texttt{PrimalDualInteriorPointProblem}):
\begin{itemize}
\item the norm of the infeasibility measure is determined by that of the feasibility problem or that of the relaxed problem);
\item the objective measure is $\objectivemeasure{\pi}(x) \equaldef \pi f(x)$ (for all the currently implemented strategies);
\item the auxiliary measure accumulates proximal (in the feasibility problem \texttt{l1RelaxedProblem}) and barrier (in \texttt{PrimalDualInteriorPointProblem}) \\ terms.
\end{itemize}
The fact that the progress measures (respectively their local models) are closely linked to the definition of the problem (respectively of the subproblem) guarantees that, as much as possible, a given combination of strategies is meaningful.
This assumes that i) the approximate solution to the subproblem satisfies a sufficient decrease condition (e.g., it achieves a fraction of the Cauchy decrease) or that the direction is a descent direction for the KKT residuals, and ii) that the reformulations are (asymptotically) equivalent to the original problem.
}

\subsubsection{\revision{Inertia correction}}
\label{sec:uno-inertia-correction}

\revision{
Each Hessian model (see Section~\ref{sec:hessian-models}) declares whether it is positive definite (the exact Hessian is ``indefinite'' by default), and whether it is available as an explicit matrix and as a linear operator. The inertia of the Lagrangian Hessian is corrected if:
\begin{itemize}
\item the Hessian model is not already positive definite; and
\item the Hessian model is available as an explicit matrix (otherwise, Uno terminates with an error message); and
\item the inertia correction strategy performs primal correction (``primal'' or ``primal-dual'').
\end{itemize}
}

\revision{
Furthermore, the successive reformulations (\texttt{l1RelaxedProblem} and \\\texttt{PrimalDualInteriorPointProblem}) maintain the lists of variables and constraints that should be primally and dually regularized, respectively. The expected inertia is computed accordingly.
For instance, \texttt{l1RelaxedProblem} adds elastic variables that enter the objective and constraints only linearly. Consequently, they contribute to the number of zero eigenvalues.
}

\subsubsection{\revision{Instantiation of ingredients}}
\label{sec:instantation-ingredients}

\revision{
Most ingredients are picked by the user via options. The following ingredients may be set or overridden by Uno after analyzing the 
\begin{itemize}
\item for inequality-constrained methods, the subproblem solver is picked depending whether the subproblems have curvature. If this is the case, a QP solver is allocated. Otherwise, an LP solver is allocated;
\item if the reformulated problem is unconstrained,
    \begin{itemize}
    \item the constraint relaxation strategy is set to the custom \texttt{NoRelaxation} strategy;
    \item the globalization strategy is set to \texttt{l1MeritFunction};
    \item and there is no curvature in the subproblem, the subproblem solver is set to the custom \texttt{BoxLPSolver};
    \end{itemize}
\item if the Hessian model is ``exact'' but neither an explicit Hessian matrix nor a Hessian linear operator was provided by the modeler, the Hessian model is set to ``zero'' (as long as quasi-Newton approximations are not implemented in Uno) and a warning message is printed.
\end{itemize}
}

\subsubsection{\revision{Sparse matrix formats}}

\revision{
Each subproblem solver possesses an object that inherits from the abstract class \texttt{EvaluationSpace} in which they store the Jacobian matrix, the Hessian matrix, or the augmented matrix in specific sparse formats:
\begin{itemize}
\item BQPD expects a Jacobian in Compressed Sparse Row (CSR) format;
\item MA57, MA27, and MUMPS expect a matrix in COOrdinate (COO) format;
\item HiGHS expects the Jacobian in Compressed Sparse (Row or Column) format.
\end{itemize}
To avoid redundancy, these matrices exist in a single sparse format and at a single location in the code. The \texttt{EvaluationSpace} object is passed around throughout the framework and provides Hessian-vector, Jacobian-vector, and Jacobian-transposed-vector products (used for instance for computing the stationarity residual).
}

\subsubsection{Termination criteria}

In practice, we cannot hope to drive the error in the first-order conditions to zero: optimization algorithms may terminate at locally infeasible points or at points that fail to satisfy constraint qualifications. We therefore check for termination, rather than for optimality.
Uno terminates at the primal-dual iterate $(x^*, \constraintmultipliers^*, \boundmultipliers^*, \pi^*)$ if
\begin{itemize}
\item sufficient first-order optimality conditions are approximately satisfied:
   \begin{itemize}
   \item a \textbf{feasible KKT point} (CQ holds) if it satisfies \eqref{eq:stationarity}, \eqref{eq:primal-feasibility}, \eqref{eq:dual-feasibility}, and \eqref{eq:complementarity} with $\pi^* > 0$;
   \item a \textbf{feasible FJ point} (CQ does not hold) if it satisfies \eqref{eq:stationarity}, \eqref{eq:primal-feasibility}, \eqref{eq:dual-feasibility}, and \eqref{eq:complementarity} with $\pi^* = 0$;
   \item an \textbf{infeasible stationary point} (a minimum of the constraint violation) if it satisfies \eqref{eq:stationarity} and \eqref{eq:dual-feasibility} with $\pi^* = 0$, no primal feasibility $c(x^*) > 0$ and a complementarity condition on the violated constraints that depends on the norm used in the feasibility problem;
   \end{itemize}
\item primal feasibility is approximately satisfied and the trust-region radius is close to machine epsilon. \revision{This may occur when the problem is non-differentiable, when the gradients are erroneous, or when second-order derivatives are not used, as is the case in Sequential Linear Programming (SLP) methods};
\item the first-order optimality conditions cannot be satisfied for the user-defined tolerance $\varepsilon$ but are satisfied for a looser tolerance (e.g., $100 \varepsilon$) for a certain number of consecutive iterations (e.g., 15).
\end{itemize}

\subsubsection{Error handling}

Uno \revision{terminates with an error message} if it encounters an IEEE exception at the initial point $x^{(0)}$. Otherwise, it tries to recover from IEEE exceptions during the optimization process by invoking the globalization mechanism in the following way:
\begin{itemize}
\item the current trust-region radius is reduced;
\item the backtracking line search reduces the tentative step length.
\end{itemize}



%% file: sections/combining_the_ingredients.tex
\section{Combining the ingredients}
\label{sec:combining-the-ingredients}

We show below how the \revision{eight} ingredients naturally arise in a popular \revision{trust-region restoration filter SQP method} and a line-search restoration filter interior-point method.
\revision{
We also provide complete multicolor pseudocodes in Algorithms~\ref{alg:uno-tr-filter-sqp} and \ref{alg:uno-ls-filter-interior-point}.
They highlight the complexity of algorithm design, and the contrast between the basic (simplistic) abstract algorithm (Algorithm~\ref{alg:iterative-methods}) and a robust and convergent optimization algorithm. In addition, the images illustrate the challenges we faced when implementing Uno.
}

\subsection{Trust-region restoration filter SQP}

Convergence of SQP filter methods has been proven under mild conditions in \cite{Fletcher2002Global,FleLeyToi:02} in the context of trust-region methods and in \cite{WaeBie:05a,WaeBie:05b} in the context of line-search methods.
The trust-region optimality QP subproblem about $(x^{(k)}, \constraintmultipliers^{(k)})$ is defined as
\begin{equation}
\tag{$\mathit{QP}^{(k)}(\Delta^{(l)})$}
\begin{array}{lll} \dps
\mini_{\highlight[subproblem_color]{\direction{x}}} & \highlight[subproblem_color]{\frac{1}{2} \direction{x}^T H_{\highlight[constraint_relaxation_color]{1}}^{(k)} \direction{x} + \highlight[constraint_relaxation_color]{1} (\nabla f^{(k)})^T \direction{x}} \\
\st 	& \highlight[subproblem_color]{c^{(k)} + (\nabla c^{(k)})^T \direction{x} = 0} \\
		& \highlight[subproblem_color]{x^{(k)} + \direction{x} \ge 0} \\
		& \highlight[mechanism_color]{\lVert \direction{x} \rVert_\infty \le \Delta^{(l)}},
\end{array}
\label{eq:TR-QP}
\end{equation}
where $\Delta^{(l)} > 0$ is the current trust-region radius.
If \eqref{eq:TR-QP} is infeasible (\revision{$\Delta^{(l)}$ is too small or} the linearized constraints are inconsistent), we switch to feasibility restoration and solve a smooth reformulation of the $\ell_1$ feasibility problem with elastic variables $\positiveelastic \in \mathbb{R}^\numberconstraints$ and $\negativeelastic \in \mathbb{R}^\numberconstraints$: 
\begin{equation}
\tag{$\mathit{FQP}^{(k)}(\Delta^{(l)})$}
\begin{array}{lll} \dps
\mini_{\highlight[subproblem_color]{\direction{x}}, \highlight[constraint_relaxation_color]{\positiveelastic,\negativeelastic}} & \highlight[subproblem_color]{\frac{1}{2} \direction{x}^T H_{\highlight[constraint_relaxation_color]{0}}^{(k)} \direction{x}} ~\highlight[constraint_relaxation_color]{+~ e^T \positiveelastic + e^T \negativeelastic} \\
\st 	& \highlight[subproblem_color]{c^{(k)} + (\nabla c^{(k)})^T \direction{x} ~\highlight[constraint_relaxation_color]{-~ \positiveelastic + \negativeelastic} = 0} \\
        & \highlight[subproblem_color]{x^{(k)} + \direction{x} \ge 0} \\
        & \highlight[mechanism_color]{\lVert \direction{x} \rVert_\infty \le \Delta^{(l)}} \\
        & \highlight[constraint_relaxation_color]{\positiveelastic \ge 0, ~ \negativeelastic \ge 0}.
\end{array}
\end{equation}

If the trial iterate $x^{(k)} + \direction{x}$ makes sufficient progress with respect to the filter method, it is accepted and we start a new iteration; if \revision{the trust region} was active at the solution of the QP ($\| \direction{x}^* \|_{\infty} = \Delta^{(l)}$), we enlarge the \revision{radius}.
If the trial iterate is rejected, we resolve the trust-region subproblem with a smaller trust-region radius. It can be shown that this mechanism either generates an acceptable iterate or results in an infeasible QP.
The complete pseudocode of the method is given in Algorithm~\ref{alg:uno-tr-filter-sqp}.

Note that in \cite{Fletcher2002Global} the authors do not solve the feasibility problem with elastic variables but exploit the partition into satisfied and violated linearized constraints provided by the Phase I method of the QP solver BQPD.

\subsection{Line-search filter restoration \revision{infeasible} interior-point method}

\revision{
An infeasible interior-point method does not require feasibility with respect to the general constraints. A prerequisite is to turn inequality constraints into equality constraints using slack variables:
\begin{equation*}
\begin{array}{ll} \dps
\mini_x  & f(x) \\
\st 	 & l_c \le c(x) \le u_c \\
         & l_x \le x \le u_x
\end{array}
\qquad \Longleftrightarrow \qquad
\begin{array}{ll} \dps
\mini_{x,s}  & f(x) \\
\st 	 & c(x) - s = 0 \\
         & l_x \le x \le u_x \\
         & l_c \le s \le u_c,
\end{array}
\end{equation*}
before applying the reformulation described in Section~\ref{sec:ipm}.
}

Provided that the subproblem is convex, the primal-dual direction is the solution of a Newton linear system, the \textit{primal-dual system} (or full-space system), in which the multipliers $\boundmultipliers$ are treated as independent variables. We solve a smaller, symmetrized version of \revision{the primal-dual system}:
\begin{equation}
\tag{$IPSP_{\barrierparameter}$}
\begin{aligned}
\begin{pmatrix}
\highlight[subproblem_color]{H_{\highlight[constraint_relaxation_color]{1}}^{(k)} + (X^{(k)})^{-1} \Boundmultipliers^{(k)}} \highlight[mechanism_color]{+~ \delta_w I}~   & \highlight[subproblem_color]{\nabla c^{(k)}} \\
\highlight[subproblem_color]{(\nabla c^{(k)})^T}                              & \highlight[mechanism_color]{-\delta_c I} \\
\end{pmatrix}
\begin{pmatrix}
\direction{x} \\
-\direction{\constraintmultipliers}
\end{pmatrix}
= \\
-
\begin{pmatrix}
\highlight[subproblem_color]{\highlight[constraint_relaxation_color]{1}\nabla f^{(k)} - \nabla c^{(k)} \constraintmultipliers^{(k)} - \barrierparameter (X^{(k)})^{-1} e} \\
\highlight[subproblem_color]{c^{(k)}}
\end{pmatrix},
\end{aligned}
\end{equation}
where $X^{(k)} = \text{diag}(x^{(k)})$, $\Boundmultipliers^{(k)} = \text{diag}(\boundmultipliers^{(k)})$, $e$ is a vector of ones of appropriate size, $\delta_w$ and $\delta_c$ are primal and dual \revision{inertia correction} coefficients, respectively. The dual direction for the bound constraints is given by $\direction{\boundmultipliers} = (X^{(k)})^{-1} (\barrierparameter e - \Boundmultipliers^{(k)} \direction{x}) - \boundmultipliers^{(k)}$.
The fraction-to-boundary rule determines primal and dual step lengths that maintain positivity of $x$ and $\boundmultipliers$:
\begin{equation}
\label{eq:fraction-to-boundary}
\begin{aligned}
\alpha_x^{(k)} \equaldef \max \{ \alpha \in (0, 1] ~|~ x^{(k)} + \alpha \direction{x} \ge (1 - \tau) x^{(k)} \} \\
\alpha_{\boundmultipliers}^{(k)} \equaldef \max \{ \alpha \in (0, 1] ~|~ \boundmultipliers^{(k)} + \alpha \direction{\boundmultipliers} \ge (1 - \tau) \boundmultipliers^{(k)} \},
\end{aligned}
\end{equation}
where $\tau$ is a parameter close to 1. 
A filter line search
assesses whether the trial iterate makes sufficient progress with respect to the filter method. If so, the trial iterate is accepted; otherwise the line search backtracks and tries again with a smaller step length.
If the step length ultimately falls below a given threshold (e.g., $10^{-7}$), we switch to feasibility restoration and solve the $\ell_1$ feasibility problem with elastic variables.


Note that by construction the filter entries depend on the barrier parameter $\barrierparameter$ through the auxiliary measure $\auxiliarymeasure$. Consequently, the filter must be flushed whenever $\barrierparameter$ is updated.
The complete pseudocode of the method is given in Algorithm~\ref{alg:uno-ls-filter-interior-point}.
An alternative that we plan to explore in the future is to update the filter \revision{entries} whenever the barrier parameter is updated.

%% file: sections/results.tex
\section{Numerical results}
\label{sec:results}

We now compare the \revision{\texttt{filtersqp} and \texttt{ipopt} presets of Uno \unoversion} against state-of-the-art solvers \revision{filterSQP (20010817, with the QP solver BQPD), IPOPT 3.14.20 (with the linear solver MA57 3.10.0), SNOPT 7.5-1.2, MINOS 5.51, LANCELOT, LOQO 7.03, and CONOPT~3.17A} on 429 small test problems of the CUTE benchmark~\cite{bongartz1995cute} translated into AMPL\footnote{The AMPL translations are available at \url{https://arnold-neumaier.at/glopt/coconut/Benchmark/Library2_new_v1.html}}. \revision{Their dimensions are reported in Table~\ref{tab:cute-dimensions} (the full name of the \texttt{nuffield} is  \texttt{nuffield\_continuum}, which the table is too narrow to contain).}
The log files of all the solvers are available in the following repository: \\
\revision{\url{https://github.com/cvanaret/nonlinear_optimization_solver_benchmark}}~.
IPOPT terminated with the status ``EXIT: Problem has too few degrees of freedom'' on the instances \texttt{argauss} and \texttt{lewispol}; this is interpreted as failure.

\renewcommand{\arraystretch}{0.85}
\setlength{\tabcolsep}{3pt}
\begin{table}[htbp]
\centering
\caption{\revision{Dimensions of the 429 small CUTE instances.}}
\resizebox{\textwidth}{!}{%
\begin{tabular}{|l|c|c||l|c|c||l|c|c||l|c|c||l|c|c|}
\hline
instance & $n$ & $m$ & instance & $n$ & $m$ & instance & $n$ & $m$ & instance & $n$ & $m$ & instance & $n$ & $m$ \\
\hline
3pk & 30 & 0 & dual1 & 85 & 1 & hs021 & 2 & 1 & hs110 & 10 & 0 & palmer3 & 4 & 0 \\
aircrfta & 5 & 5 & dual2 & 96 & 1 & hs022 & 2 & 2 & hs111lnp & 10 & 3 & palmer4a & 6 & 0 \\
aircrftb & 5 & 0 & dual4 & 75 & 1 & hs023 & 2 & 5 & hs111 & 10 & 3 & palmer4b & 4 & 0 \\
airport & 84 & 42 & dualc1 & 9 & 13 & hs024 & 2 & 2 & hs112 & 10 & 3 & palmer4c & 8 & 0 \\
aljazzaf & 3 & 1 & dualc2 & 7 & 9 & hs025 & 3 & 0 & hs113 & 10 & 8 & palmer4e & 8 & 0 \\
allinitc & 3 & 1 & dualc5 & 8 & 1 & hs026 & 3 & 1 & hs114 & 10 & 11 & palmer4 & 4 & 0 \\
allinit & 3 & 0 & dualc8 & 8 & 15 & hs027 & 3 & 1 & hs116 & 13 & 15 & palmer5a & 8 & 0 \\
allinitu & 4 & 0 & eg1 & 3 & 0 & hs028 & 3 & 1 & hs117 & 15 & 5 & palmer5b & 9 & 0 \\
alsotame & 2 & 1 & eigencco & 30 & 15 & hs029 & 3 & 1 & hs118 & 15 & 17 & palmer5c & 6 & 0 \\
argauss & 3 & 15 & eigmaxc & 22 & 22 & hs030 & 3 & 1 & hs119 & 16 & 8 & palmer5d & 4 & 0 \\
arglinb & 10 & 0 & eigminc & 22 & 22 & hs031 & 3 & 1 & hs21mod & 7 & 1 & palmer5e & 8 & 0 \\
arglinc & 8 & 0 & engval2 & 3 & 0 & hs032 & 3 & 2 & hs268 & 5 & 5 & palmer6a & 6 & 0 \\
avgasa & 8 & 10 & errinros & 50 & 0 & hs033 & 3 & 2 & hs35mod & 2 & 1 & palmer6c & 8 & 0 \\
avgasb & 8 & 10 & expfita & 5 & 21 & hs034 & 3 & 2 & hs3mod & 2 & 0 & palmer6e & 8 & 0 \\
avion2 & 49 & 15 & expfit & 2 & 0 & hs035 & 3 & 1 & hs44new & 4 & 5 & palmer7a & 6 & 0 \\
bard & 3 & 0 & extrasim & 2 & 1 & hs036 & 3 & 1 & hs99exp & 28 & 21 & palmer7c & 8 & 0 \\
batch & 46 & 69 & extrosnb & 10 & 0 & hs037 & 3 & 1 & hubfit & 2 & 1 & palmer7e & 8 & 0 \\
beale & 2 & 0 & fccu & 19 & 8 & hs038 & 4 & 0 & humps & 2 & 0 & palmer8a & 6 & 0 \\
biggs3 & 3 & 0 & fletcher & 4 & 4 & hs039 & 4 & 2 & hypcir & 2 & 2 & palmer8c & 8 & 0 \\
biggs5 & 5 & 0 & genhs28 & 10 & 8 & hs040 & 4 & 3 & jensmp & 2 & 0 & palmer8e & 8 & 0 \\
biggs6 & 6 & 0 & genhumps & 5 & 0 & hs041 & 4 & 1 & kiwcresc & 3 & 2 & pentagon & 6 & 12 \\
biggsc4 & 4 & 7 & gigomez1 & 3 & 3 & hs042 & 3 & 1 & kowosb & 4 & 0 & pfit1ls & 3 & 0 \\
booth & 2 & 2 & goffin & 51 & 50 & hs043 & 4 & 3 & lakes & 90 & 78 & pfit1 & 3 & 0 \\
box2 & 2 & 0 & gottfr & 2 & 2 & hs044 & 4 & 6 & launch & 25 & 29 & pfit2ls & 3 & 0 \\
box3 & 3 & 0 & gridnetg & 44 & 34 & hs045 & 5 & 0 & lewispol & 6 & 9 & pfit2 & 3 & 0 \\
bqp1var & 1 & 0 & gridneth & 61 & 36 & hs046 & 5 & 2 & linspanh & 72 & 32 & pfit3ls & 3 & 0 \\
bqpgabim & 46 & 0 & gridneti & 61 & 36 & hs047 & 5 & 3 & loadbal & 31 & 31 & pfit3 & 3 & 0 \\
bqpgasim & 50 & 0 & growthls & 3 & 0 & hs048 & 5 & 2 & loghairy & 2 & 0 & pfit4ls & 3 & 0 \\
brkmcc & 2 & 0 & growth & 3 & 0 & hs049 & 5 & 2 & logros & 2 & 0 & pfit4 & 3 & 0 \\
brownal & 10 & 0 & gulf & 3 & 0 & hs050 & 5 & 3 & lootsma & 3 & 2 & polak1 & 3 & 2 \\
brownbs & 2 & 0 & hadamals & 90 & 0 & hs051 & 5 & 3 & lotschd & 12 & 7 & polak2 & 11 & 2 \\
brownden & 4 & 0 & haifas & 7 & 9 & hs052 & 5 & 3 & lsnnodoc & 5 & 4 & polak3 & 12 & 10 \\
bt10 & 2 & 2 & hairy & 2 & 0 & hs053 & 5 & 3 & lsqfit & 2 & 1 & polak4 & 3 & 3 \\
bt11 & 5 & 3 & haldmads & 6 & 42 & hs054 & 6 & 1 & madsen & 3 & 6 & polak5 & 3 & 2 \\
bt12 & 5 & 3 & hart6 & 6 & 0 & hs055 & 6 & 6 & makela1 & 3 & 2 & polak6 & 5 & 4 \\
bt13 & 5 & 1 & hatflda & 4 & 0 & hs056 & 7 & 4 & makela2 & 3 & 3 & portfl1 & 12 & 1 \\
bt1 & 2 & 1 & hatfldb & 4 & 0 & hs057 & 2 & 1 & makela3 & 21 & 20 & portfl2 & 12 & 1 \\
bt2 & 3 & 1 & hatfldc & 4 & 0 & hs059 & 2 & 3 & makela4 & 21 & 40 & portfl3 & 12 & 1 \\
bt3 & 5 & 3 & hatfldd & 3 & 0 & hs060 & 3 & 1 & maratosb & 2 & 0 & portfl4 & 12 & 1 \\
bt4 & 3 & 2 & hatflde & 3 & 0 & hs061 & 3 & 2 & maratos & 2 & 1 & portfl6 & 12 & 1 \\
bt5 & 3 & 2 & hatfldf & 3 & 3 & hs062 & 3 & 1 & matrix2 & 6 & 2 & powellbs & 2 & 2 \\
bt6 & 5 & 2 & hatfldg & 25 & 25 & hs063 & 3 & 2 & maxlika & 8 & 0 & powellsq & 2 & 2 \\
bt7 & 5 & 3 & hatfldh & 4 & 7 & hs064 & 3 & 1 & mconcon & 15 & 11 & prodpl0 & 60 & 29 \\
bt8 & 5 & 2 & heart6ls & 6 & 0 & hs065 & 3 & 1 & mdhole & 2 & 0 & prodpl1 & 60 & 29 \\
bt9 & 4 & 2 & heart6 & 6 & 6 & hs066 & 3 & 2 & methanb8 & 31 & 0 & pspdoc & 4 & 0 \\
byrdsphr & 3 & 2 & heart8ls & 8 & 0 & hs067 & 10 & 7 & methanl8 & 31 & 0 & qudlin & 12 & 0 \\
camel6 & 2 & 0 & heart8 & 8 & 8 & hs070 & 4 & 1 & mexhat & 2 & 0 & recipe & 3 & 3 \\
cantilvr & 5 & 1 & helix & 3 & 0 & hs071 & 4 & 2 & meyer3 & 3 & 0 & res & 11 & 7 \\
catena & 32 & 11 & hilberta & 10 & 0 & hs072 & 4 & 2 & mifflin1 & 3 & 2 & rk23 & 17 & 11 \\
cb2 & 3 & 3 & hilbertb & 50 & 0 & hs073 & 4 & 3 & mifflin2 & 3 & 2 & robot & 7 & 2 \\
cb3 & 3 & 3 & himmelba & 2 & 2 & hs074 & 4 & 4 & minmaxbd & 5 & 20 & rosenbr & 2 & 0 \\
chaconn1 & 3 & 3 & himmelbb & 2 & 0 & hs075 & 4 & 4 & minmaxrb & 3 & 4 & rosenmmx & 5 & 4 \\
chaconn2 & 3 & 3 & himmelbc & 2 & 2 & hs076 & 4 & 3 & minsurf & 36 & 0 & s365mod & 7 & 5 \\
chebyqad & 50 & 0 & himmelbd & 2 & 2 & hs077 & 5 & 2 & mistake & 9 & 13 & sim2bqp & 2 & 0 \\
chnrosnb & 50 & 0 & himmelbe & 3 & 3 & hs078 & 5 & 3 & model & 60 & 32 & simbqp & 2 & 0 \\
cliff & 2 & 0 & himmelbf & 4 & 0 & hs079 & 5 & 3 & mwright & 5 & 3 & simpllpa & 2 & 2 \\
cluster & 2 & 2 & himmelbg & 2 & 0 & hs080 & 5 & 3 & nasty & 2 & 0 & simpllpb & 2 & 3 \\
concon & 15 & 11 & himmelbh & 2 & 0 & hs081 & 5 & 3 & nonmsqrt & 9 & 0 & sineali & 20 & 0 \\
congigmz & 3 & 5 & himmelbk & 24 & 14 & hs083 & 5 & 3 & nuffield & 2 & 1 & sineval & 2 & 0 \\
coolhans & 9 & 9 & himmelp1 & 2 & 0 & hs084 & 5 & 3 & obstclal & 64 & 0 & sisser & 2 & 0 \\
core1 & 65 & 50 & himmelp2 & 2 & 1 & hs085 & 5 & 36 & obstclbl & 64 & 0 & snake & 2 & 2 \\
coshfun & 61 & 20 & himmelp3 & 2 & 2 & hs086 & 5 & 6 & obstclbu & 64 & 0 & spanhyd & 72 & 32 \\
cresc4 & 6 & 8 & himmelp4 & 2 & 3 & hs087 & 9 & 4 & odfits & 10 & 6 & spiral & 3 & 2 \\
csfi1 & 5 & 4 & himmelp5 & 2 & 3 & hs088 & 2 & 1 & optcntrl & 28 & 20 & ssnlbeam & 31 & 20 \\
csfi2 & 5 & 4 & himmelp6 & 2 & 4 & hs089 & 3 & 1 & optmass & 66 & 55 & stancmin & 3 & 2 \\
cube & 2 & 0 & hong & 4 & 1 & hs090 & 4 & 1 & optprloc & 30 & 29 & supersim & 2 & 2 \\
dallass & 44 & 29 & hs001 & 2 & 0 & hs091 & 5 & 1 & orthregb & 27 & 6 & swopf & 82 & 91 \\
deconvb & 51 & 0 & hs002 & 2 & 0 & hs092 & 6 & 1 & orthrege & 36 & 20 & synthes1 & 6 & 6 \\
deconvc & 51 & 1 & hs003 & 2 & 0 & hs093 & 6 & 2 & osbornea & 5 & 0 & tame & 2 & 1 \\
deconvu & 51 & 0 & hs004 & 2 & 0 & hs095 & 6 & 4 & osborneb & 11 & 0 & tointqor & 50 & 0 \\
degenlpa & 20 & 14 & hs005 & 2 & 0 & hs096 & 6 & 4 & oslbqp & 8 & 0 & try-b & 2 & 1 \\
degenlpb & 20 & 15 & hs006 & 2 & 1 & hs097 & 6 & 4 & palmer1a & 6 & 0 & twobars & 2 & 2 \\
demymalo & 3 & 3 & hs007 & 2 & 1 & hs098 & 6 & 4 & palmer1b & 4 & 0 & vanderm4 & 9 & 17 \\
denschna & 2 & 0 & hs008 & 2 & 2 & hs099 & 19 & 14 & palmer1c & 8 & 0 & watson & 31 & 0 \\
denschnb & 2 & 0 & hs009 & 2 & 1 & hs100lnp & 7 & 2 & palmer1d & 7 & 0 & weeds & 3 & 0 \\
denschnc & 2 & 0 & hs010 & 2 & 1 & hs100mod & 7 & 4 & palmer1e & 8 & 0 & womflet & 3 & 3 \\
denschnd & 3 & 0 & hs011 & 2 & 1 & hs100 & 7 & 4 & palmer1 & 4 & 0 & yfit & 3 & 0 \\
denschne & 3 & 0 & hs012 & 2 & 1 & hs101 & 7 & 6 & palmer2a & 6 & 0 & yfitu & 3 & 0 \\
denschnf & 2 & 0 & hs013 & 2 & 1 & hs102 & 7 & 6 & palmer2b & 4 & 0 & zangwil2 & 2 & 0 \\
dipigri & 7 & 4 & hs014 & 2 & 2 & hs103 & 7 & 6 & palmer2c & 8 & 0 & zangwil3 & 3 & 3 \\
disc2 & 28 & 23 & hs015 & 2 & 2 & hs104 & 8 & 6 & palmer2e & 8 & 0 & zecevic2 & 2 & 2 \\
discs & 33 & 66 & hs016 & 2 & 2 & hs105 & 8 & 0 & palmer2 & 4 & 0 & zecevic3 & 2 & 2 \\
dixchlng & 10 & 5 & hs017 & 2 & 2 & hs106 & 8 & 6 & palmer3a & 6 & 0 & zecevic4 & 2 & 2 \\
dixon3dq & 10 & 0 & hs018 & 2 & 2 & hs107 & 9 & 6 & palmer3b & 4 & 0 & zigzag & 58 & 50 \\
djtl & 2 & 0 & hs019 & 2 & 2 & hs108 & 9 & 13 & palmer3c & 8 & 0 & zy2 & 3 & 1 \\
dnieper & 57 & 24 & hs020 & 2 & 3 & hs109 & 9 & 10 & palmer3e & 8 & 0 &  & &  \\
\hline
\end{tabular}
}
\label{tab:cute-dimensions}
\end{table}

\begin{figure}[h!]
\centering
\begin{subfigure}[b]{0.49\columnwidth}
\centering
\includegraphics[width=\textwidth]{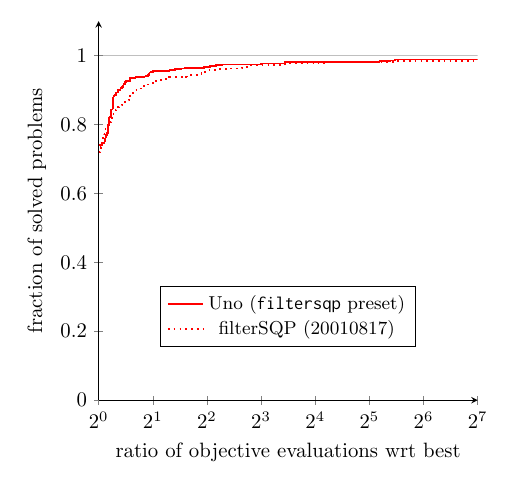}
\caption{Uno \texttt{filtersqp} vs filterSQP}
\label{fig:performance-profile-filtersqp}
\end{subfigure}
\hfill
\begin{subfigure}[b]{0.49\columnwidth}
\centering
\includegraphics[width=\textwidth]{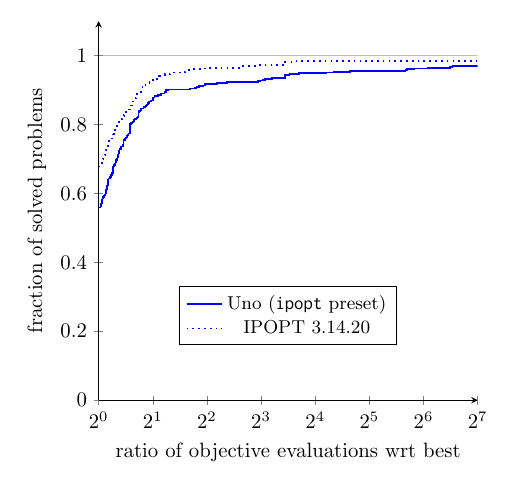}
\caption{Uno \texttt{ipopt} vs IPOPT}
\label{fig:performance-profile-ipopt}
\end{subfigure}
\begin{subfigure}[b]{0.85\columnwidth}
\centering
\includegraphics[width=\textwidth]{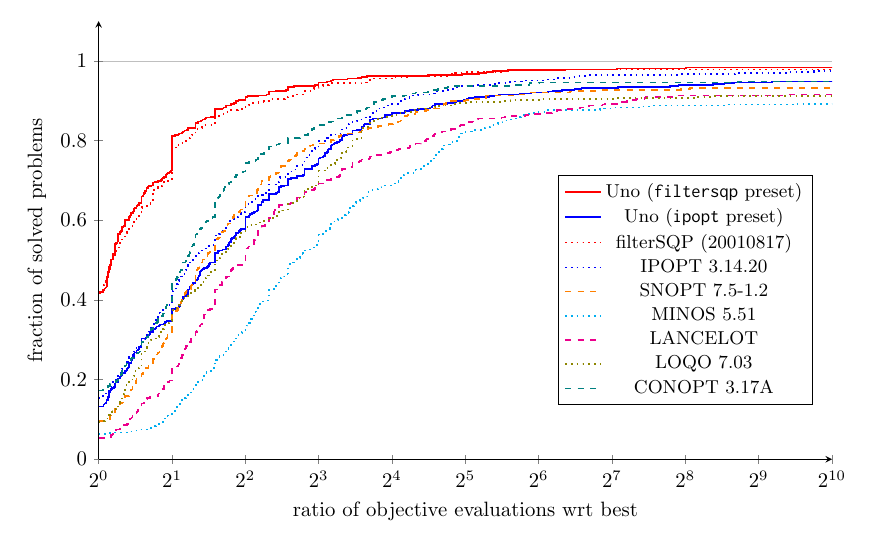}
\caption{Uno vs state-of-the-art solvers}
\label{fig:performance-profile-all}
\end{subfigure}
\caption{Performance profile on 429 small CUTE problems: number of objective evaluations.}
\label{fig:performance-profile}
\end{figure}

In this section we demonstrate that the Uno presets closely mimic the corresponding \revision{state-of-the-art} solvers \revision{filterSQP and IPOPT}. 
\revision{Figures~\ref{fig:performance-profile-filtersqp} and \ref{fig:performance-profile-ipopt} portray performance profiles~\cite{DolaMore:02} of the \texttt{filtersqp} preset vs filterSQP, and of the \texttt{ipopt} preset vs IPOPT, respectively. Figure~\ref{fig:performance-profile-all} illustrates the overall picture.}

\revision{On all three plots}, the y axis is the fraction of solved problems, and the x axis represents the relative budget of objective evaluations compared with the (virtual) best solver for each instance. The higher and more to the left, the better.
These performance profiles \revision{demonstrate that the Uno presets \texttt{filtersqp} and \texttt{ipopt} mimic the corresponding state-of-the-art solvers well and perform well with respect to all the state-of-the-art solvers.} These results validate the implementation of off-the-shelf strategies within Uno. 

Table~\ref{tab:statistics_LPs_problems} Table~\ref{tab:statistics_QPs_problems} and Table~\ref{tab:statistics_NLPs_problems} summarize the number of objective evaluations required by each solver for linear, quadratic, and nonlinear instances, respectively. The lowest count is shown in bold.
\revision{12 instances (\texttt{core1}, \texttt{csfi2}, \texttt{discs}, \texttt{growth}, \texttt{hs085}, \texttt{hs093}, \texttt{hs114}, \texttt{launch}, \texttt{meyer3}, \texttt{polak6}, \texttt{spanhyd}, \texttt{vanderm4})} are solved by IPOPT but not by the Uno \texttt{ipopt} preset, most likely because of the IPOPT features that were not implemented (see Section~\ref{sec:combinations}).

\revision{
Table~\ref{tab:geometric-mean} shows for each solver the shifted geometric mean of its objective evaluations $e_1, \ldots, e_N$: $\dps (\prod_{i=1}^N (e_i + s))^\frac{1}{N} - s$, where $s = 10$ is the shift factor, and $N = 429$ is the number of instances. Upon failure on instance $i$, we take $e_i = 10^6$.
}
\begin{table}[htbp]
\renewcommand{\arraystretch}{1.1}
\setlength{\tabcolsep}{2pt}
\caption{\revision{Shifted geometric mean ($s = 10$) of objective evaluations on 429 small CUTE instances}.}
\begin{tabular}{|cc|ccccccc|}
\hline
\multicolumn{2}{|c|}{Uno presets} & \multicolumn{7}{c|}{State-of-the-art solvers} \\
\texttt{filtersqp} & \texttt{ipopt} & filterSQP & IPOPT & SNOPT & MINOS & LANCELOT & LOQO & CONOPT \\
\hline
18.08 & 48.45 & 18.73 & 34.63 & 61.02 & 159.78 & 92.27 & 71.44 & 43.15 \\
\hline
\end{tabular}
\label{tab:geometric-mean}
\end{table}

%% file: sections/conclusion.tex
\section{Conclusion and future developments}

We have introduced an abstract framework for unifying nonlinearly constrained optimization with \revision{eight} key ingredients common to most methods. We have shown that our abstract framework provides common notation and abstractions for the unified description of well-known optimization strategies. 
We then presented Uno, a C\texttt{++} implementation of the \revision{unifying} framework. A modular solver, Uno provides efficient implementations of off-the-shelf strategies and facilitates the development of new algorithmic ideas by making building blocks and abstractions readily available. In particular, newly developed strategies (e.g., a new \revision{globalization strategy}) can be immediately deployed and tested within a wide range of strategy combinations.
\revision{Strategy combinations that do not exist as software implementations} can be generated on the fly with no programming effort and tested on a given instance.
We believe that Uno has the potential to serve as an experimentation laboratory for the optimization community and accelerate research in \revision{nonlinearly constrained} optimization.

Future releases of Uno will include the following features: quasi-Newton methods (L-BFGS and \revision{SR1}), iterative linear solvers (e.g., MINRES~\cite{paige1975solution}), equality-constrained methods (SLP-EQP), \revision{SQuID's S$\ell_1$QP  method}~\cite{burke2014sequential}, parallel line search~\cite{schittkowski2001nlpqlp}, and additional interfaces to subproblem solvers.




%% file: sections/appendices.tex
\section*{Appendices}

\input{sections/statistics_table}

\input{sections/solvers}
\input{sections/uml}
\input{sections/combining_filtersqp}
\input{sections/combining_ipopt}

%% file: sections/statistics_table.tex
\begin{scriptsize}
\renewcommand{\arraystretch}{1.1}
\setlength{\tabcolsep}{2pt}
\begin{longtable}{|c|cc|ccccccc|}
\caption{Number of objective evaluations on a subset of \revision{CUTE} LPs.}
\label{tab:statistics_LPs_problems} \\
\hline
 & \multicolumn{2}{c|}{Uno presets}           & \multicolumn{7}{c|}{State-of-the-art solvers} \\
Problem & \texttt{filtersqp} & \texttt{ipopt} & filterSQP & IPOPT & SNOPT & MINOS & LANCELOT & LOQO & CONOPT \\
\hline
\endfirsthead
\hline
 & \multicolumn{2}{c|}{Uno presets}           & \multicolumn{7}{c|}{State-of-the-art solvers} \\
Problem & \texttt{filtersqp} & \texttt{ipopt} & filterSQP & IPOPT & SNOPT & MINOS & LANCELOT & LOQO & CONOPT \\
\hline
\endhead
booth & 2 & 2 & 2 & 2 & 2 & \textbf{1} & 3 & 7 & \textbf{1} \\
degenlpa & 2 & 44 & 2 & 29 & 26 & 15 & 25 & 29 & \textbf{1} \\
degenlpb & 2 & 42 & 2 & 41 & 26 & 23 & 45 & 30 & \textbf{1} \\
extrasim & 2 & 6 & 2 & 6 & \textbf{1} & \textbf{1} & 3 & 11 & \textbf{1} \\
goffin & 2 & 7 & 3 & 8 & 25 & 25 & 10 & 13 & \textbf{1} \\
himmelba & 2 & 2 & 2 & 2 & 2 & \textbf{1} & 3 & 7 & \textbf{1} \\
linspanh & 2 & 18 & 2 & 19 & 5 & 14 & 13 & 12 & \textbf{1} \\
makela4 & 3 & 8 & 3 & 8 & \textbf{1} & \textbf{1} & 20 & 12 & \textbf{1} \\
model & 2 & 22 & 2 & 15 & 23 & 34 & 33 & 15 & \textbf{1} \\
res & 2 & 9 & 2 & 10 & \textbf{1} & 5 & \textbf{1} & 12 & \textbf{1} \\
simpllpa & 2 & 13 & 2 & 14 & 3 & 3 & 5 & 12 & \textbf{1} \\
simpllpb & 2 & 10 & 2 & 11 & \textbf{1} & \textbf{1} & 4 & 13 & \textbf{1} \\
supersim & 2 & 7 & 2 & 2 & \textbf{1} & \textbf{1} & 7 & 11 & \textbf{1} \\
zangwil3 & 2 & 2 & 2 & 2 & 3 & 2 & 3 & 8 & \textbf{1} \\
\hline
\end{longtable}
\end{scriptsize}

\begin{scriptsize}
\renewcommand{\arraystretch}{1.1}
\setlength{\tabcolsep}{2pt}
\begin{longtable}{|c|cc|ccccccc|}
\caption{Number of objective evaluations on a subset of \revision{CUTE} QPs.}
\label{tab:statistics_QPs_problems} \\
\hline
 & \multicolumn{2}{c|}{Uno presets}           & \multicolumn{7}{c|}{State-of-the-art solvers} \\
Problem & \texttt{filtersqp} & \texttt{ipopt} & filterSQP & IPOPT & SNOPT & MINOS & LANCELOT & LOQO & CONOPT \\
\hline
\endfirsthead
\hline
 & \multicolumn{2}{c|}{Uno presets}           & \multicolumn{7}{c|}{State-of-the-art solvers} \\
Problem & \texttt{filtersqp} & \texttt{ipopt} & filterSQP & IPOPT & SNOPT & MINOS & LANCELOT & LOQO & CONOPT \\
\hline
\endhead
3pk & \textbf{7} & 11 & \textbf{7} & 12 & 41 & 209 & 47 & 20 & 24 \\
arglinb & 2 & 3 & 2 & 3 & \textbf{1} & 8 & 2 & 3 & 20 \\
arglinc & 2 & 5 & 2 & 3 & \textbf{1} & 8 & 2 & 3 & 18 \\
avgasa & \textbf{2} & 10 & \textbf{2} & 10 & 9 & 19 & 12 & 11 & 16 \\
avgasb & \textbf{2} & 13 & \textbf{2} & 13 & 9 & 15 & 11 & 16 & 19 \\
biggsc4 & \textbf{2} & 26 & \textbf{2} & 36 & 13 & 16 & 15 & 21 & 4 \\
bqp1var & 2 & 6 & 2 & 6 & \textbf{1} & 5 & 3 & 10 & 2 \\
bqpgabim & \textbf{2} & 17 & \textbf{2} & 21 & 36 & 107 & 6 & 14 & 17 \\
bqpgasim & \textbf{2} & 17 & \textbf{2} & 21 & 40 & 121 & 6 & 14 & 15 \\
bt3 & \textbf{2} & \textbf{2} & \textbf{2} & \textbf{2} & 5 & 12 & 6 & 3 & 9 \\
deconvb & 42 & 8514 & 30 & -- & 103 & 30 & 31 & 45 & \textbf{23} \\
dixon3dq & \textbf{2} & \textbf{2} & \textbf{2} & \textbf{2} & 10 & 43 & 4 & 3 & 9 \\
dual1 & \textbf{2} & 17 & \textbf{2} & 17 & 221 & -- & 9 & 22 & 29 \\
dual2 & \textbf{2} & 14 & \textbf{2} & 14 & 116 & -- & 9 & 19 & 27 \\
dual4 & \textbf{2} & 14 & \textbf{2} & 14 & 31 & -- & 8 & 21 & 23 \\
dualc1 & \textbf{2} & 30 & \textbf{2} & 29 & 7 & 24 & 17 & 26 & 14 \\
dualc2 & \textbf{2} & 21 & \textbf{2} & 27 & 5 & 14 & 19 & 21 & 17 \\
dualc5 & \textbf{2} & 11 & \textbf{2} & 11 & 7 & 24 & 8 & 16 & 15 \\
dualc8 & \textbf{2} & 15 & \textbf{2} & 13 & 9 & 16 & 24 & 19 & 25 \\
fccu & 4 & \textbf{2} & 4 & \textbf{2} & 19 & 48 & 13 & 3 & 11 \\
genhs28 & \textbf{2} & \textbf{2} & \textbf{2} & \textbf{2} & 11 & 12 & 5 & 3 & 9 \\
hatfldc & \textbf{5} & 6 & \textbf{5} & 6 & 14 & 19 & 9 & 9 & 15 \\
hatfldh & \textbf{2} & 14 & \textbf{2} & 19 & 3 & 8 & 13 & 16 & 5 \\
hilberta & \textbf{2} & \textbf{2} & \textbf{2} & \textbf{2} & \textbf{2} & 12 & 60 & 3 & 9 \\
hilbertb & \textbf{2} & \textbf{2} & \textbf{2} & \textbf{2} & 50 & -- & 5 & 3 & 10 \\
hs003 & \textbf{2} & 6 & \textbf{2} & 5 & \textbf{2} & 12 & 16 & 11 & 5 \\
hs021 & 2 & 7 & 2 & 9 & \textbf{1} & 8 & 3 & 12 & 12 \\
hs028 & \textbf{2} & \textbf{2} & \textbf{2} & \textbf{2} & 4 & 12 & 4 & 3 & 8 \\
hs035 & \textbf{2} & 8 & \textbf{2} & 8 & 5 & 12 & 8 & 10 & 12 \\
hs044 & \textbf{2} & 20 & \textbf{2} & 21 & \textbf{2} & 10 & 11 & 16 & 5 \\
hs048 & \textbf{2} & \textbf{2} & \textbf{2} & \textbf{2} & 6 & 16 & 3 & 3 & 9 \\
hs051 & \textbf{2} & \textbf{2} & \textbf{2} & \textbf{2} & 6 & 12 & 10 & 3 & 7 \\
hs052 & \textbf{2} & \textbf{2} & \textbf{2} & \textbf{2} & 5 & 11 & 6 & 3 & 9 \\
hs053 & \textbf{2} & 7 & \textbf{2} & 7 & \textbf{2} & 12 & 6 & 12 & 15 \\
hs054 & \textbf{2} & 8 & \textbf{2} & 8 & 5 & 16 & 9 & 12 & 12 \\
hs076 & \textbf{2} & 8 & \textbf{2} & 8 & 4 & 13 & 9 & 11 & 16 \\
hs118 & \textbf{3} & 12 & \textbf{3} & 12 & 21 & 41 & 19 & 17 & 4 \\
hs21mod & 2 & 16 & 2 & 17 & \textbf{1} & 8 & 3 & 20 & 12 \\
hs268 & \textbf{2} & 19 & \textbf{2} & 17 & 6 & 36 & 27 & 27 & 19 \\
hs35mod & 2 & 16 & 2 & 16 & \textbf{1} & 6 & 3 & 16 & 6 \\
hs3mod & \textbf{2} & 6 & \textbf{2} & 6 & 5 & 13 & 4 & 12 & 10 \\
hs44new & \textbf{2} & 14 & \textbf{2} & 14 & 4 & 10 & 10 & 21 & 4 \\
lotschd & \textbf{2} & 15 & 3 & 15 & 8 & 5 & 9 & 19 & 7 \\
lsqfit & \textbf{2} & 8 & \textbf{2} & 8 & 3 & 11 & 7 & 12 & 10 \\
maratosb & 10 & 33 & 13 & 33 & \textbf{7} & 10 & 8 & 21 & 14 \\
nasty & 2 & 2 & 2 & 2 & \textbf{1} & 6 & 5 & 3 & 10 \\
obstclal & \textbf{2} & 11 & \textbf{2} & 15 & 37 & 90 & 7 & 23 & 15 \\
obstclbl & \textbf{2} & 12 & \textbf{2} & 13 & 44 & 87 & 3 & 13 & 47 \\
obstclbu & \textbf{2} & 13 & \textbf{2} & 13 & 36 & 74 & \textbf{2} & 15 & 26 \\
oslbqp & \textbf{2} & 15 & \textbf{2} & 15 & 6 & 11 & 3 & 19 & \textbf{2} \\
palmer1c & 7 & 23 & 7 & \textbf{2} & 8 & 61 & 145 & -- & -- \\
palmer1d & 6 & 26 & 6 & \textbf{2} & 7 & 50 & 34 & -- & -- \\
palmer2c & 5 & 3 & 5 & \textbf{2} & 8 & 62 & 298 & -- & -- \\
palmer3c & 5 & \textbf{2} & 5 & \textbf{2} & 8 & 61 & 206 & -- & -- \\
palmer4c & 6 & \textbf{2} & 6 & \textbf{2} & 8 & 61 & 176 & -- & -- \\
palmer5c & 4 & \textbf{2} & 4 & \textbf{2} & 6 & 26 & \textbf{2} & -- & -- \\
palmer5d & 5 & \textbf{2} & 5 & \textbf{2} & 4 & 25 & \textbf{2} & -- & -- \\
palmer6c & 6 & \textbf{2} & 6 & \textbf{2} & 8 & 60 & 159 & -- & -- \\
palmer7c & 8 & 3 & 8 & \textbf{2} & 8 & 64 & 189 & -- & -- \\
palmer8c & 7 & \textbf{2} & 7 & \textbf{2} & 8 & 61 & 152 & -- & -- \\
portfl1 & \textbf{2} & 10 & \textbf{2} & 10 & 12 & 48 & 20 & 17 & 20 \\
portfl2 & \textbf{2} & 9 & \textbf{2} & 9 & 12 & 50 & 14 & 17 & 23 \\
portfl3 & \textbf{2} & 11 & \textbf{2} & 11 & 13 & 52 & 14 & 17 & 32 \\
portfl4 & \textbf{2} & 10 & \textbf{2} & 10 & 11 & 49 & 19 & 16 & 19 \\
portfl6 & \textbf{2} & 9 & \textbf{2} & 9 & 11 & 48 & 20 & 16 & 31 \\
qudlin & \textbf{2} & 24 & \textbf{2} & 26 & 11 & 18 & \textbf{2} & 20 & 3 \\
sim2bqp & \textbf{2} & 8 & \textbf{2} & 8 & \textbf{2} & 9 & 3 & 14 & 5 \\
simbqp & \textbf{2} & 8 & \textbf{2} & 8 & \textbf{2} & 10 & \textbf{2} & 13 & 4 \\
tame & 2 & 6 & 2 & 6 & \textbf{1} & 8 & 2 & 9 & 6 \\
tointqor & \textbf{2} & \textbf{2} & \textbf{2} & \textbf{2} & 50 & -- & 8 & 3 & 12 \\
zangwil2 & \textbf{2} & \textbf{2} & \textbf{2} & \textbf{2} & \textbf{2} & 11 & 4 & 3 & 3 \\
zecevic2 & \textbf{2} & 9 & \textbf{2} & 9 & \textbf{2} & 9 & 7 & 11 & 7 \\
\hline
\end{longtable}
\end{scriptsize}

\begin{scriptsize}
\renewcommand{\arraystretch}{1.1}
\setlength{\tabcolsep}{2pt}
\begin{longtable}{|c|cc|ccccccc|}
\caption{Number of objective evaluations on a subset of \revision{CUTE} NLPs.}
\label{tab:statistics_NLPs_problems} \\
\hline
 & \multicolumn{2}{c|}{Uno presets}           & \multicolumn{7}{c|}{State-of-the-art solvers} \\
Problem & \texttt{filtersqp} & \texttt{ipopt} & filterSQP & IPOPT & SNOPT & MINOS & LANCELOT & LOQO & CONOPT \\
\hline
\endfirsthead
\hline
 & \multicolumn{2}{c|}{Uno presets}           & \multicolumn{7}{c|}{State-of-the-art solvers} \\
Problem & \texttt{filtersqp} & \texttt{ipopt} & filterSQP & IPOPT & SNOPT & MINOS & LANCELOT & LOQO & CONOPT \\
\hline
\endhead
aircrfta & 4 & 4 & 4 & 4 & -- & \textbf{1} & 10 & 6 & 2 \\
aircrftb & 21 & \textbf{15} & 21 & 19 & 58 & 65 & 27 & 21 & 32 \\
airport & 14 & 21 & \textbf{13} & 16 & 58 & 527 & 69 & 19 & 101 \\
aljazzaf & 11 & 79 & 15 & 82 & 145 & 64 & 24 & 79 & \textbf{9} \\
allinitc & 29 & 34 & 24 & 44 & 105 & 55 & 76 & -- & \textbf{17} \\
allinit & \textbf{11} & 15 & \textbf{11} & 19 & 17 & 29 & 13 & 18 & 16 \\
allinitu & 12 & 15 & 12 & 15 & 14 & 20 & 15 & \textbf{11} & 17 \\
alsotame & 5 & 9 & 5 & 9 & 6 & 11 & 11 & 11 & \textbf{4} \\
argauss & \textbf{1} & -- & \textbf{1} & -- & 7 & -- & -- & -- & 3 \\
avion2 & 23 & 6918 & 19 & 136 & 19 & \textbf{17} & 787 & 65 & 25 \\
bard & 11 & \textbf{9} & 11 & \textbf{9} & 23 & 36 & 15 & 19 & 17 \\
batch & \textbf{9} & 62 & \textbf{9} & 34 & 33 & 379 & -- & 64 & 32 \\
beale & 19 & 11 & \textbf{10} & 19 & 15 & 24 & 21 & \textbf{10} & 16 \\
biggs3 & \textbf{11} & 16 & \textbf{11} & 28 & 24 & 31 & 40 & 12 & 15 \\
biggs5 & 54 & \textbf{24} & 50 & 36 & 107 & 38 & 64 & 32 & 48 \\
biggs6 & 74 & \textbf{38} & 83 & 50 & 120 & 119 & 103 & 59 & 139 \\
box2 & \textbf{9} & \textbf{9} & \textbf{9} & \textbf{9} & 10 & 10 & 20 & 10 & 10 \\
box3 & \textbf{8} & 11 & \textbf{8} & 15 & 24 & 15 & 31 & 11 & 13 \\
brkmcc & \textbf{4} & \textbf{4} & \textbf{4} & \textbf{4} & 10 & 13 & 7 & 8 & 10 \\
brownal & \textbf{8} & \textbf{8} & \textbf{8} & \textbf{8} & 21 & 77 & 24 & 10 & 17 \\
brownbs & 48 & 9 & 48 & 8 & 32 & 35 & 7 & 49 & \textbf{1} \\
brownden & \textbf{9} & \textbf{9} & \textbf{9} & \textbf{9} & 40 & 41 & \textbf{9} & 57 & 15 \\
bt10 & 8 & 9 & 7 & 7 & \textbf{1} & \textbf{1} & 20 & 11 & 3 \\
bt11 & 9 & 10 & \textbf{7} & 9 & 12 & 56 & 22 & 11 & 21 \\
bt12 & \textbf{5} & \textbf{5} & \textbf{5} & \textbf{5} & 10 & 59 & 11 & 11 & 11 \\
bt13 & 24 & 26 & 48 & 25 & 34 & 69 & -- & \textbf{21} & 590 \\
bt1 & 2 & 15 & \textbf{1} & 15 & 12 & 16 & 19 & 24 & \textbf{1} \\
bt2 & \textbf{13} & 14 & \textbf{13} & \textbf{13} & 18 & 385 & 36 & 18 & 20 \\
bt4 & \textbf{9} & 12 & 11 & 10 & 10 & 47 & 25 & \textbf{9} & 21 \\
bt5 & 10 & 9 & 9 & \textbf{8} & 11 & 172 & 20 & \textbf{8} & 10 \\
bt6 & \textbf{11} & 12 & 12 & 18 & 14 & 117 & 25 & 12 & 28 \\
bt7 & 17 & 117 & 19 & 30 & 36 & 85 & 49 & 18 & \textbf{8} \\
bt8 & 23 & 31 & 12 & 50 & 14 & 21 & 30 & -- & \textbf{9} \\
bt9 & 19 & 31 & 23 & \textbf{14} & 19 & 34 & 22 & \textbf{14} & 17 \\
byrdsphr & \textbf{8} & 98 & 11 & 19 & 59 & 35 & 43 & 11 & 17 \\
camel6 & 9 & 11 & \textbf{8} & 11 & 19 & 24 & \textbf{8} & 12 & 19 \\
cantilvr & 14 & 16 & 16 & \textbf{12} & 23 & 69 & 27 & 16 & 45 \\
catena & 13 & 54 & 13 & \textbf{7} & 87 & 125 & 56 & 26 & 80 \\
cb2 & 8 & 16 & \textbf{7} & 9 & \textbf{7} & 31 & 18 & 11 & 14 \\
cb3 & 8 & 15 & 7 & 10 & \textbf{1} & 26 & 18 & 11 & 3 \\
chaconn1 & 6 & 13 & \textbf{5} & 7 & 9 & 22 & 12 & 11 & 10 \\
chaconn2 & 6 & 12 & 5 & 7 & \textbf{1} & 20 & 11 & 11 & 4 \\
chebyqad & 64 & 111 & 50 & 139 & -- & -- & 62 & \textbf{12} & 87 \\
chnrosnb & 56 & \textbf{55} & 59 & 92 & 170 & -- & 68 & 58 & 101 \\
cliff & 28 & 28 & 28 & \textbf{24} & 28 & 46 & 28 & 32 & 32 \\
cluster & 10 & 10 & 10 & 10 & -- & \textbf{3} & 45 & 11 & \textbf{3} \\
concon & \textbf{3} & 10 & 5 & 10 & 14 & -- & 676 & -- & 40 \\
congigmz & 5 & 40 & 4 & 33 & 10 & 16 & 30 & 36 & \textbf{3} \\
coolhans & 4 & 10 & 3 & 10 & -- & 4 & 281 & 25 & \textbf{2} \\
core1 & \textbf{5} & -- & 6 & 86 & 44 & 86 & -- & 52 & \textbf{5} \\
coshfun & 86 & 41226 & 303 & 725 & 273 & 1091 & 154 & \textbf{24} & 140 \\
cresc4 & \textbf{16} & 5311 & 52 & 79 & 93 & 864 & -- & 111 & 47 \\
csfi1 & 35 & 28 & 18 & 12 & 36 & -- & 155 & 16 & \textbf{7} \\
csfi2 & 11 & -- & 8 & 86 & 60 & -- & 180 & 17 & \textbf{3} \\
cube & 39 & 38 & 41 & 58 & 42 & 66 & 52 & 41 & \textbf{36} \\
dallass & \textbf{18} & 36 & 56 & 29 & 109 & 140 & -- & -- & -- \\
deconvc & \textbf{24} & 194 & 58 & 99 & 81 & 87 & 43 & 31 & 30 \\
deconvu & 4033 & 94 & 971 & 729 & 152 & \textbf{30} & 69 & 76 & 86 \\
demymalo & 9 & 19 & \textbf{8} & 12 & \textbf{8} & 34 & 28 & 15 & 15 \\
denschna & \textbf{7} & \textbf{7} & \textbf{7} & \textbf{7} & 12 & 23 & 13 & 9 & 10 \\
denschnb & 10 & 21 & 10 & 25 & 10 & 17 & 11 & \textbf{9} & 10 \\
denschnc & \textbf{11} & \textbf{11} & \textbf{11} & \textbf{11} & 21 & 30 & 13 & 14 & 13 \\
denschnd & 41 & 45 & 43 & \textbf{27} & 77 & 111 & 65 & 44 & 64 \\
denschne & \textbf{11} & 12 & \textbf{11} & 25 & 44 & 34 & 16 & 15 & 34 \\
denschnf & \textbf{7} & \textbf{7} & \textbf{7} & \textbf{7} & 12 & 21 & 8 & 11 & 12 \\
dipigri & 13 & 37 & 13 & 22 & 23 & 129 & 63 & \textbf{12} & 30 \\
disc2 & 4 & 219 & 25 & 48 & 800 & -- & -- & 28 & \textbf{3} \\
discs & 10 & -- & 40 & 269 & -- & -- & 422 & 59 & \textbf{3} \\
dixchlng & \textbf{11} & 12 & 12 & \textbf{11} & 31 & -- & 44 & 26 & 43 \\
djtl & -- & 70 & 29 & 861 & -- & -- & 100 & 132 & \textbf{1} \\
dnieper & \textbf{4} & 35 & \textbf{4} & 31 & 13 & 36 & 75 & 25 & 6 \\
eg1 & \textbf{8} & \textbf{8} & \textbf{8} & \textbf{8} & 9 & 14 & 9 & 10 & 13 \\
eigencco & 20 & \textbf{13} & 29 & 14 & 34 & 158 & 17 & 22 & 26 \\
eigmaxc & \textbf{1} & 9 & 7 & 7 & 18 & -- & 21 & 14 & 3 \\
eigminc & \textbf{1} & 11 & 7 & 8 & 21 & 182 & 11 & 13 & 4 \\
engval2 & \textbf{19} & 21 & 20 & 33 & 34 & 66 & 30 & 28 & 33 \\
errinros & 53 & \textbf{49} & 53 & 70 & 267 & -- & 76 & 60 & 90 \\
expfita & \textbf{13} & 34 & \textbf{13} & 31 & 24 & 28 & 54 & 23 & 14 \\
expfit & 13 & \textbf{9} & 13 & \textbf{9} & 18 & 28 & 11 & 13 & 19 \\
extrosnb & 2 & 2 & 2 & \textbf{1} & 2 & 5 & \textbf{1} & \textbf{1} & \textbf{1} \\
fletcher & \textbf{1} & 29 & \textbf{1} & 28 & 2 & -- & 28 & 14 & 19 \\
genhumps & 139 & 257 & 188 & 234 & 77 & 169 & 134 & 146 & \textbf{71} \\
gigomez1 & \textbf{8} & 19 & \textbf{8} & 19 & 9 & 51 & 33 & 16 & 14 \\
gottfr & 6 & 8 & 13 & 9 & -- & \textbf{1} & 35 & 12 & 3 \\
gridnetg & 5 & 17 & \textbf{4} & 11 & 22 & 42 & 21 & 12 & 21 \\
gridneth & \textbf{5} & 7 & \textbf{5} & 7 & 36 & 114 & 20 & 12 & 12 \\
gridneti & \textbf{5} & 18 & \textbf{5} & 14 & 47 & 111 & 22 & 16 & 24 \\
growthls & \textbf{101} & \textbf{101} & 106 & 170 & 184 & 258 & 178 & 122 & 212 \\
growth & \textbf{101} & -- & 106 & 157 & 187 & 258 & 178 & 143 & 206 \\
gulf & \textbf{26} & 28 & \textbf{26} & 44 & 66 & 695 & 63 & \textbf{26} & 32 \\
hadamals & \textbf{13} & 104 & \textbf{13} & 138 & 19 & 266 & 20 & \textbf{13} & 38 \\
haifas & 12 & 36 & 13 & \textbf{10} & 28 & 59 & 27 & 13 & 23 \\
hairy & 45 & 51 & 84 & 115 & \textbf{35} & 59 & 102 & 64 & 92 \\
haldmads & \textbf{11} & 103 & 41 & 259 & 71 & 94 & 45 & 31 & 12 \\
hart6 & 15 & 10 & 11 & 14 & 16 & 36 & \textbf{9} & 23 & 20 \\
hatflda & 10 & 11 & 15 & 11 & 30 & 44 & 47 & \textbf{9} & 114 \\
hatfldb & \textbf{9} & 11 & 11 & 11 & 28 & 32 & 28 & 11 & 138 \\
hatfldd & 31 & \textbf{23} & \textbf{23} & 27 & 29 & 48 & 66 & 25 & 37 \\
hatflde & 36 & 28 & \textbf{26} & 32 & 31 & 63 & 57 & 29 & 31 \\
hatfldf & 4 & 3151 & 15 & 1335 & -- & \textbf{2} & 113 & 13 & 3 \\
hatfldg & 4 & 20 & 15 & 20 & -- & \textbf{2} & 29 & 16 & 3 \\
heart6ls & 2596 & \textbf{985} & -- & 1544 & -- & -- & -- & -- & 3067 \\
heart6 & \textbf{1} & 1516 & 16 & 75 & -- & 49 & -- & 327 & 3 \\
heart8ls & 248 & \textbf{108} & 217 & 189 & -- & 474 & 238 & \textbf{108} & 369 \\
heart8 & \textbf{3} & 37 & 12 & 40 & -- & -- & 359 & 39 & \textbf{3} \\
helix & 18 & \textbf{13} & 19 & 25 & 28 & 53 & 18 & \textbf{13} & 21 \\
himmelbb & 62 & 15 & 25 & 12 & \textbf{8} & 20 & 11 & 16 & 32 \\
himmelbc & 5 & 8 & 8 & 9 & -- & \textbf{1} & 11 & 8 & 3 \\
himmelbd & 4 & 7327 & 4 & 79 & 12 & -- & -- & -- & \textbf{3} \\
himmelbe & 3 & 3 & 2 & 3 & -- & \textbf{1} & 8 & 8 & 2 \\
himmelbf & \textbf{8} & 12 & \textbf{8} & 11 & 53 & 47 & 29 & 24 & 15 \\
himmelbg & 10 & 10 & 10 & 14 & 12 & 17 & 16 & \textbf{7} & 13 \\
himmelbh & 8 & 20 & 8 & 24 & 10 & 10 & \textbf{7} & 9 & 10 \\
himmelbk & \textbf{6} & 21 & \textbf{6} & 19 & 82 & 109 & 206 & 23 & 11 \\
himmelp1 & \textbf{9} & 12 & \textbf{9} & 12 & 19 & 20 & 29 & 11 & 21 \\
himmelp2 & \textbf{9} & 22 & \textbf{9} & 19 & 32 & 151 & 275 & 20 & 18 \\
himmelp3 & \textbf{5} & 10 & \textbf{5} & 13 & 8 & 119 & 870 & 17 & 9 \\
himmelp4 & \textbf{5} & 10 & \textbf{5} & 25 & 8 & 115 & 737 & 17 & 9 \\
himmelp5 & \textbf{12} & 25 & \textbf{12} & 158 & 44 & 75 & 273 & 90 & 23 \\
himmelp6 & 2 & 11 & 2 & 12 & 2 & 5 & 2 & 45 & \textbf{1} \\
hong & 8 & 16 & 5 & 13 & \textbf{4} & 14 & 6 & 20 & \textbf{4} \\
hs001 & 33 & 34 & 36 & 53 & 48 & \textbf{9} & 41 & 35 & 36 \\
hs002 & 9 & 14 & 9 & 17 & 15 & 12 & \textbf{7} & 21 & 10 \\
hs004 & 3 & 6 & 3 & 6 & 4 & 6 & \textbf{2} & 8 & \textbf{2} \\
hs005 & \textbf{8} & 9 & 11 & 9 & 9 & 13 & 9 & 10 & 14 \\
hs006 & \textbf{3} & \textbf{3} & \textbf{3} & 7 & 9 & 90 & 63 & 10 & 4 \\
hs007 & 13 & 65 & 13 & 28 & 30 & 64 & 26 & \textbf{12} & 13 \\
hs008 & 6 & 6 & 6 & 6 & -- & \textbf{1} & 13 & 8 & 3 \\
hs009 & \textbf{5} & \textbf{5} & \textbf{5} & 6 & 10 & 11 & 22 & 9 & 8 \\
hs010 & \textbf{10} & 26 & \textbf{10} & 13 & 20 & 22 & 18 & 15 & 19 \\
hs011 & 7 & 17 & \textbf{6} & 9 & 15 & 46 & 16 & 12 & 12 \\
hs012 & 9 & 16 & \textbf{8} & 9 & 11 & 158 & 26 & 11 & 16 \\
hs013 & 26 & 45 & 34 & 78 & \textbf{17} & 53 & 60 & -- & 19 \\
hs014 & 7 & 12 & 6 & 8 & 10 & 9 & 13 & 11 & \textbf{3} \\
hs015 & \textbf{3} & 25 & 7 & 21 & 11 & 85 & 47 & 31 & \textbf{3} \\
hs016 & 6 & 12 & 5 & 23 & 5 & 9 & 19 & 18 & \textbf{3} \\
hs017 & 9 & 27 & 8 & 18 & 19 & 11 & 20 & 30 & \textbf{6} \\
hs018 & 8 & 18 & \textbf{7} & 27 & 32 & 93 & 117 & 15 & 16 \\
hs019 & \textbf{7} & 18 & \textbf{7} & 16 & 9 & 56 & 45 & 18 & 9 \\
hs020 & 5 & 14 & 5 & 7 & 5 & 8 & 23 & 24 & \textbf{3} \\
hs022 & 6 & 11 & \textbf{2} & 7 & 7 & 47 & 10 & 9 & 3 \\
hs023 & 7 & 15 & 7 & 12 & 7 & 53 & 51 & 18 & \textbf{6} \\
hs024 & \textbf{3} & 11 & \textbf{3} & 12 & 8 & 8 & 14 & 13 & \textbf{3} \\
hs025 & 31 & 40 & 27 & 44 & 2 & 5 & \textbf{1} & 20 & \textbf{1} \\
hs026 & 18 & 25 & 18 & 26 & 27 & 76 & 41 & \textbf{14} & 33 \\
hs027 & 23 & 70 & \textbf{8} & 151 & 21 & 136 & 31 & 16 & 43 \\
hs029 & 9 & 15 & \textbf{8} & 9 & 14 & 173 & 18 & 10 & 13 \\
hs030 & \textbf{2} & 26 & \textbf{2} & 26 & 5 & 29 & 8 & 9 & 20 \\
hs031 & 7 & 13 & \textbf{6} & 8 & 11 & 28 & 12 & 17 & 20 \\
hs032 & \textbf{2} & 16 & \textbf{2} & 20 & 5 & 14 & 7 & 24 & 3 \\
hs033 & 6 & 14 & 5 & 16 & 9 & 38 & 9 & 11 & \textbf{3} \\
hs034 & 9 & 13 & 8 & 10 & \textbf{5} & 9 & 21 & 15 & 16 \\
hs036 & \textbf{3} & 13 & \textbf{3} & 13 & 10 & 8 & 7 & 20 & \textbf{3} \\
hs037 & \textbf{6} & 13 & \textbf{6} & 13 & 10 & 16 & 13 & 11 & 15 \\
hs038 & 53 & 58 & 54 & 78 & 101 & 88 & 56 & \textbf{13} & 94 \\
hs039 & 19 & 31 & 23 & \textbf{14} & 19 & 34 & 22 & \textbf{14} & 17 \\
hs040 & 6 & 6 & 5 & \textbf{4} & 9 & 21 & 11 & 8 & 11 \\
hs041 & 3 & 9 & \textbf{2} & 11 & 7 & 5 & 7 & 16 & 14 \\
hs042 & 7 & 10 & \textbf{6} & 7 & 10 & 21 & 13 & 9 & 13 \\
hs043 & \textbf{10} & 15 & 11 & \textbf{10} & 11 & 102 & 25 & 12 & 27 \\
hs045 & 2 & 24 & 2 & 48 & 2 & 5 & \textbf{1} & 25 & \textbf{1} \\
hs046 & 19 & 24 & 19 & 20 & 32 & 121 & 28 & \textbf{18} & 37 \\
hs047 & \textbf{21} & 22 & \textbf{21} & \textbf{21} & 28 & 125 & 29 & 24 & 31 \\
hs049 & \textbf{17} & 21 & \textbf{17} & 20 & 34 & 61 & 38 & 21 & 21 \\
hs050 & \textbf{9} & 10 & \textbf{9} & 10 & 20 & 24 & 11 & 16 & 14 \\
hs055 & 2 & 8 & 2 & 8 & 2 & 5 & 7 & 11 & \textbf{1} \\
hs056 & \textbf{3} & 18 & 19 & 39 & 52 & 52 & 12 & 12 & 27 \\
hs057 & 5 & 24 & 5 & 28 & -- & 28 & \textbf{2} & 14 & 28 \\
hs059 & \textbf{10} & 91 & 11 & 72 & 22 & 127 & 340 & 27 & 22 \\
hs060 & 8 & \textbf{7} & \textbf{7} & 8 & 13 & 120 & 18 & 9 & 23 \\
hs061 & \textbf{1} & 11 & \textbf{1} & 10 & 39 & 70 & 19 & 10 & 2 \\
hs062 & \textbf{8} & 10 & 10 & 9 & 16 & 18 & 37 & 13 & 19 \\
hs063 & 8 & 11 & \textbf{1} & 8 & 18 & 119 & 21 & 8 & 19 \\
hs064 & \textbf{11} & 19 & 13 & 18 & 28 & 98 & 38 & 26 & 20 \\
hs065 & 6 & 15 & \textbf{5} & 91 & 11 & 319 & 42 & 19 & 22 \\
hs066 & 9 & 11 & 14 & \textbf{8} & \textbf{8} & \textbf{8} & 11 & 15 & 15 \\
hs067 & 13 & 18 & \textbf{12} & \textbf{12} & 32 & 60 & 287 & 17 & 19 \\
hs070 & 40 & \textbf{18} & 42 & 36 & 34 & 66 & 39 & 23 & 31 \\
hs071 & 7 & 13 & \textbf{6} & 9 & 8 & 55 & 16 & 13 & 14 \\
hs072 & \textbf{14} & 18 & 15 & 17 & 27 & 52 & 65 & 28 & 20 \\
hs073 & \textbf{4} & 12 & \textbf{4} & 9 & 11 & 8 & 18 & 21 & 5 \\
hs074 & \textbf{6} & 15 & \textbf{6} & 10 & 15 & 27 & 13 & 16 & 15 \\
hs075 & \textbf{5} & 15 & \textbf{5} & 10 & 12 & 18 & 110 & 18 & 6 \\
hs077 & \textbf{12} & \textbf{12} & 14 & 13 & 16 & 123 & 27 & \textbf{12} & 28 \\
hs078 & 6 & 6 & \textbf{5} & \textbf{5} & 7 & 57 & 12 & 8 & 13 \\
hs079 & 6 & 6 & \textbf{5} & \textbf{5} & 12 & 54 & 14 & 7 & 24 \\
hs080 & 9 & 9 & 8 & \textbf{7} & 9 & 48 & 13 & 10 & 17 \\
hs081 & 30 & 12 & 38 & \textbf{8} & 11 & 55 & 17 & 16 & 17 \\
hs083 & 6 & 19 & \textbf{5} & 15 & 8 & 11 & 16 & 14 & 6 \\
hs084 & 11 & 23 & 6 & 12 & 58 & 45 & 47 & 44 & \textbf{4} \\
hs085 & -- & -- & -- & 143 & -- & -- & -- & 63 & \textbf{7} \\
hs086 & \textbf{5} & 13 & \textbf{5} & 11 & 16 & 13 & 15 & 13 & 19 \\
hs087 & \textbf{7} & 33 & \textbf{7} & 18 & 16 & 23 & 32 & 25 & 11 \\
hs088 & 26 & 44 & 19 & 18 & 59 & 66 & 56 & 27 & \textbf{17} \\
hs089 & 21 & 21 & 31 & 38 & 85 & 197 & 61 & 30 & \textbf{19} \\
hs090 & 73 & 46 & \textbf{2} & 28 & 55 & 92 & 58 & 29 & 54 \\
hs091 & 110 & 122 & 337 & 15 & 73 & 215 & 62 & 28 & \textbf{12} \\
hs092 & 83 & 41 & \textbf{2} & 25 & 56 & 110 & 58 & 22 & 103 \\
hs093 & 3 & -- & \textbf{2} & 10 & 33 & -- & -- & 13 & 22 \\
hs095 & 2 & 27 & 3 & 18 & \textbf{1} & \textbf{1} & 24 & 17 & 4 \\
hs096 & 2 & 24 & 3 & 24 & \textbf{1} & \textbf{1} & 23 & 22 & 4 \\
hs097 & \textbf{7} & 24 & \textbf{7} & 24 & 13 & 63 & 19 & 19 & 9 \\
hs098 & \textbf{7} & 43 & \textbf{7} & 21 & 13 & 47 & 19 & 93 & 9 \\
hs099 & -- & 12 & 9 & \textbf{7} & 19 & 60 & -- & 22 & 29 \\
hs100lnp & 13 & 26 & 14 & 21 & 33 & 133 & 32 & \textbf{12} & 36 \\
hs100mod & \textbf{13} & 30 & 14 & 27 & 32 & 123 & 137 & 15 & 32 \\
hs100 & 13 & 37 & 13 & 22 & 23 & 129 & 63 & \textbf{12} & 30 \\
hs101 & \textbf{16} & 375 & 34 & 273 & 530 & -- & -- & 64 & 49 \\
hs102 & \textbf{17} & 60 & 42 & 36 & 238 & 980 & -- & 154 & 46 \\
hs103 & 33 & 76 & \textbf{28} & 64 & 177 & 1418 & -- & 88 & 38 \\
hs104 & 19 & 15 & 23 & 11 & 29 & 85 & -- & 14 & \textbf{3} \\
hs105 & \textbf{9} & 21 & \textbf{9} & 31 & 89 & 114 & -- & 17 & 24 \\
hs106 & 16 & 17 & 17 & \textbf{15} & 34 & -- & -- & 27 & 33 \\
hs107 & 7 & 17 & \textbf{6} & 12 & 14 & 20 & 26 & 35 & 21 \\
hs108 & 35 & 29 & 36 & \textbf{17} & 152 & 164 & 43 & 20 & 32 \\
hs109 & \textbf{6} & 60 & 7 & 44 & 349 & -- & -- & 45 & 13 \\
hs110 & \textbf{5} & 7 & \textbf{5} & 7 & 11 & 43 & \textbf{5} & 8 & 15 \\
hs111lnp & 41 & 25 & 31 & \textbf{16} & 64 & 388 & 57 & 17 & 30 \\
hs111 & 41 & 25 & 31 & 16 & 70 & 388 & 46 & \textbf{15} & 29 \\
hs112 & 13 & 18 & \textbf{12} & 18 & 35 & 92 & 47 & 19 & 58 \\
hs113 & \textbf{6} & 16 & \textbf{6} & 12 & 28 & 146 & 97 & 17 & 30 \\
hs114 & \textbf{1} & -- & \textbf{1} & 75 & 9 & -- & 664 & -- & 4 \\
hs116 & \textbf{12} & 29 & 14 & 26 & 75 & 52 & -- & 24 & 25 \\
hs117 & 7 & 25 & \textbf{6} & 23 & 20 & 157 & 66 & 19 & 46 \\
hs119 & 8 & 17 & \textbf{7} & 15 & 22 & 29 & 28 & 29 & 18 \\
hs99exp & 13 & 117 & 12 & 30 & 42 & 212 & -- & 256 & \textbf{4} \\
hubfit & \textbf{2} & 9 & \textbf{2} & 9 & 8 & 11 & 8 & -- & -- \\
humps & \textbf{114} & 321 & -- & 249 & 257 & 193 & -- & 307 & 259 \\
hypcir & 5 & 6 & 8 & 8 & -- & \textbf{1} & 10 & 8 & 3 \\
jensmp & 11 & 11 & 11 & \textbf{10} & 36 & 55 & \textbf{10} & 14 & 13 \\
kiwcresc & 12 & 26 & \textbf{11} & \textbf{11} & 13 & 30 & 23 & 14 & 17 \\
kowosb & 18 & 15 & 18 & 23 & 33 & 39 & 24 & \textbf{11} & 26 \\
lakes & 12 & 18 & 63 & 19 & 39 & -- & -- & 288 & \textbf{1} \\
launch & \textbf{1} & -- & \textbf{1} & 409 & 246 & -- & -- & -- & 3 \\
lewispol & \textbf{1} & -- & \textbf{1} & -- & 6 & -- & -- & -- & \textbf{1} \\
loadbal & \textbf{8} & 17 & \textbf{8} & 18 & 58 & 130 & 62 & 23 & 26 \\
loghairy & 90 & 2071 & -- & -- & 249 & 402 & -- & \textbf{64} & 150 \\
logros & \textbf{50} & 95 & \textbf{50} & 358 & 109 & 147 & 66 & 398 & 84 \\
lootsma & 6 & 14 & 5 & 16 & 9 & 38 & 9 & 12 & \textbf{3} \\
lsnnodoc & 7 & 16 & 7 & 15 & 8 & 7 & 11 & 21 & \textbf{3} \\
madsen & \textbf{14} & 26 & 25 & 25 & \textbf{14} & 28 & 26 & 26 & 28 \\
makela1 & 14 & 20 & 15 & 19 & \textbf{8} & 27 & 19 & 15 & 20 \\
makela2 & \textbf{5} & 11 & \textbf{5} & 8 & 16 & 21 & 40 & 12 & 7 \\
makela3 & 22 & 35 & 25 & 17 & 288 & 96 & 125 & 18 & \textbf{15} \\
maratos & 11 & 42 & 10 & \textbf{5} & 9 & 20 & 9 & 7 & 11 \\
matrix2 & 23 & 45 & \textbf{12} & 21 & 14 & 65 & 13 & 26 & 161 \\
maxlika & \textbf{9} & 21 & \textbf{9} & 31 & 89 & 114 & -- & 17 & 24 \\
mconcon & \textbf{3} & 10 & 5 & 10 & 14 & -- & 676 & -- & 40 \\
mdhole & \textbf{3} & 70 & 56 & 106 & 70 & 116 & 68 & 98 & 6 \\
methanb8 & 47 & \textbf{8} & 47 & 9 & 306 & 175 & 221 & 155 & 28 \\
methanl8 & 105 & \textbf{66} & 96 & 82 & 499 & 669 & 640 & 80 & 84 \\
mexhat & 10 & 5 & 10 & 5 & 37 & 18 & \textbf{4} & 7 & 8 \\
meyer3 & -- & -- & \textbf{280} & 492 & -- & 775 & 559 & 522 & 716 \\
mifflin1 & 11 & 24 & 23 & \textbf{7} & 10 & 15 & 18 & 9 & 8 \\
mifflin2 & 11 & 30 & \textbf{10} & 16 & 14 & 53 & 50 & 13 & 13 \\
minmaxbd & 44 & 88 & \textbf{9} & 78 & 115 & 160 & 592 & 31 & 38 \\
minmaxrb & \textbf{3} & 11 & \textbf{3} & 11 & 4 & 34 & 81 & 14 & 5 \\
minsurf & 14 & 13 & \textbf{10} & 21 & 24 & 196 & 15 & 13 & 15 \\
mistake & 22 & 30 & 18 & \textbf{15} & 19 & 159 & 30 & \textbf{15} & 36 \\
mwright & \textbf{10} & \textbf{10} & 12 & 11 & 12 & 37 & 19 & 11 & 19 \\
nonmsqrt & -- & -- & 716 & -- & -- & 1490 & \textbf{167} & -- & -- \\
nuffield & \textbf{5} & 7 & \textbf{5} & 7 & 9 & 15 & 15 & 10 & 17 \\
odfits & \textbf{5} & 11 & 7 & 11 & 17 & 28 & 49 & 15 & 16 \\
optcntrl & \textbf{2} & 192 & 4 & 134 & 5 & 15 & 352 & 58 & 6 \\
optmass & 9 & 26 & 18 & 23 & \textbf{2} & 331 & -- & 15 & 27 \\
optprloc & 16 & 23 & \textbf{6} & 19 & 12 & 685 & 438 & 23 & 19 \\
orthregb & \textbf{2} & 3 & \textbf{2} & 3 & 9 & 262 & 64 & 8 & \textbf{2} \\
orthrege & 7995 & 41 & 180 & 166 & \textbf{31} & 857 & 795 & 482 & 13038 \\
osbornea & 4110 & \textbf{50} & -- & 152 & 120 & 127 & 57 & -- & 52 \\
osborneb & 20 & 21 & 20 & 25 & 82 & 127 & 45 & \textbf{19} & 42 \\
palmer1a & \textbf{42} & 59 & 51 & 71 & 205 & -- & 102 & 113 & 200 \\
palmer1b & \textbf{21} & \textbf{21} & \textbf{21} & 26 & 87 & -- & 55 & 65 & 69 \\
palmer1e & 105 & \textbf{49} & 74 & 153 & 186 & 149 & 353 & 187 & 100 \\
palmer1 & 33 & 1129 & 33 & 1844 & 30 & 39 & 28 & 41 & \textbf{21} \\
palmer2a & \textbf{68} & 189 & \textbf{68} & 261 & 115 & 197 & 211 & 169 & 83 \\
palmer2b & \textbf{16} & 21 & \textbf{16} & 34 & 61 & -- & 77 & 52 & 49 \\
palmer2e & 99 & \textbf{50} & 86 & 52 & 191 & 345 & 133 & 171 & 111 \\
palmer2 & 32 & 3257 & 33 & 63 & 44 & -- & 34 & \textbf{24} & 57 \\
palmer3a & \textbf{78} & 759 & 82 & 200 & 136 & -- & 201 & 172 & 172 \\
palmer3b & 21 & 18 & 21 & \textbf{15} & 54 & -- & 36 & 37 & 49 \\
palmer3e & 449 & \textbf{112} & 117 & 116 & 293 & 426 & -- & 234 & 118 \\
palmer3 & \textbf{11} & 348 & 12 & 537 & 13 & \textbf{11} & 58 & 32 & 23 \\
palmer4a & 52 & 174 & 52 & 133 & 109 & -- & 93 & 110 & \textbf{42} \\
palmer4b & 21 & \textbf{19} & 21 & 31 & 52 & -- & 64 & 35 & 48 \\
palmer4e & 25 & 39 & \textbf{24} & 38 & 123 & 210 & 123 & 87 & 65 \\
palmer4 & 12 & 363 & 12 & 1092 & 14 & \textbf{11} & 142 & 35 & 32 \\
palmer5a & \textbf{7809} & 10195 & -- & -- & -- & 321078 & -- & -- & -- \\
palmer5b & 837 & \textbf{97} & 855 & 171 & -- & 3728 & 961 & -- & -- \\
palmer5e & \textbf{3} & 5877 & \textbf{3} & -- & -- & 25501 & 8 & -- & -- \\
palmer6a & \textbf{132} & 378 & 137 & 275 & 202 & 270 & 277 & -- & -- \\
palmer6e & \textbf{22} & 39 & 38 & 59 & 198 & 259 & 47 & -- & -- \\
palmer7a & -- & \textbf{5811} & -- & -- & -- & -- & -- & -- & -- \\
palmer7e & 1526 & 9348 & -- & -- & -- & 924 & \textbf{39} & -- & -- \\
palmer8a & \textbf{49} & 105 & 50 & 102 & 127 & 103 & 61 & -- & -- \\
palmer8e & \textbf{29} & 31 & \textbf{29} & 31 & 92 & 121 & 86 & -- & -- \\
pentagon & 15 & 19 & \textbf{12} & 19 & 15 & 22 & 48 & 38 & 23 \\
pfit1ls & 365 & 395 & 566 & 682 & 480 & 717 & 453 & \textbf{337} & 1064 \\
pfit1 & 365 & 395 & 566 & 682 & 480 & 717 & 453 & \textbf{337} & 1064 \\
pfit2ls & 178 & 168 & 210 & 172 & 175 & 237 & 217 & \textbf{115} & 2586 \\
pfit2 & 178 & 168 & 210 & 172 & 175 & 237 & 217 & \textbf{115} & 2586 \\
pfit3ls & 166 & 162 & \textbf{139} & 347 & 289 & 474 & 272 & 148 & 3428 \\
pfit3 & 166 & 162 & \textbf{139} & 347 & 289 & 474 & 272 & 148 & 3428 \\
pfit4ls & \textbf{84} & 286 & 126 & 532 & 470 & 754 & 410 & 286 & 4925 \\
pfit4 & \textbf{84} & 286 & 126 & 532 & 470 & 754 & 410 & 286 & 4925 \\
polak1 & 9 & 24451 & 8 & \textbf{7} & 15 & 30 & 37 & 14 & 24 \\
polak2 & 45 & 119 & 10 & 15 & 106 & 138 & 321 & 24 & \textbf{2} \\
polak3 & 21 & -- & 24 & -- & 183 & -- & 234 & 24 & \textbf{2} \\
polak4 & \textbf{5} & 175 & \textbf{5} & 7 & 6 & 22 & 16 & 11 & 37 \\
polak5 & 81 & 55 & 45 & 33 & 43 & 16 & \textbf{7} & 69 & 17 \\
polak6 & \textbf{27} & -- & 29 & 660 & 62 & 108 & 644 & \textbf{27} & 40 \\
powellbs & 9 & 221 & 16 & 12 & -- & \textbf{1} & -- & 17 & 3 \\
powellsq & -- & 160 & 4 & 109 & 78 & -- & 23 & -- & \textbf{2} \\
prodpl0 & 10 & 20 & \textbf{9} & 16 & 60 & 42 & 35 & 24 & 34 \\
prodpl1 & 8 & 22 & 7 & 17 & 59 & 56 & 27 & 21 & \textbf{5} \\
pspdoc & 9 & 11 & \textbf{7} & 15 & 15 & 25 & 10 & 12 & 12 \\
recipe & 3 & 3 & 2 & 3 & -- & \textbf{1} & 12 & 10 & 2 \\
rk23 & 7 & 14 & 9 & 12 & 18 & 28 & 54 & 12 & \textbf{6} \\
robot & \textbf{1} & 19 & 45 & 10 & 18 & 377 & 33 & 18 & 64 \\
rosenbr & 29 & 29 & 29 & 45 & 45 & \textbf{9} & 36 & 26 & 29 \\
rosenmmx & 31 & 71 & 37 & 22 & 41 & 66 & 226 & \textbf{17} & 75 \\
s365mod & \textbf{19} & 38 & 86 & 43 & 31 & 539 & -- & 28 & 57 \\
sineali & \textbf{10} & 5792 & -- & -- & -- & 7333 & -- & 16 & 24 \\
sineval & \textbf{57} & 66 & 62 & 110 & 94 & 123 & 75 & \textbf{57} & 77 \\
sisser & 19 & 21 & 19 & 21 & \textbf{15} & 17 & 38 & 17 & 19 \\
snake & 3 & 10 & 3 & 14 & \textbf{2} & 9 & -- & -- & 4 \\
spanhyd & \textbf{6} & -- & 11 & 26 & 13 & -- & 26 & -- & 24 \\
spiral & 121 & 5743 & 152 & \textbf{64} & 130 & 148 & 96 & 135 & 642 \\
ssnlbeam & 6 & 28 & \textbf{5} & 22 & 36 & 124 & 39 & 60 & 44 \\
stancmin & \textbf{2} & 11 & \textbf{2} & 11 & 5 & 5 & 9 & 20 & 4 \\
swopf & \textbf{6} & 24 & \textbf{6} & 17 & 160 & 100 & 290 & 21 & 22 \\
synthes1 & 6 & 12 & \textbf{5} & 10 & 10 & 21 & 13 & 17 & 13 \\
try-b & 9 & 11 & 8 & 20 & 10 & 10 & 13 & 16 & \textbf{3} \\
twobars & 9 & 13 & \textbf{8} & 10 & 11 & 31 & 13 & 10 & 14 \\
vanderm4 & \textbf{1} & -- & \textbf{1} & 248 & 86 & -- & 36 & -- & 3 \\
watson & 22 & \textbf{14} & 21 & \textbf{14} & 172 & 117 & 44 & 18 & 18 \\
weeds & 38 & 26 & 39 & 32 & 51 & -- & \textbf{3} & 29 & 44 \\
womflet & 17 & 110 & \textbf{9} & 12 & 14 & 40 & 97 & 12 & 25 \\
yfit & 47 & 60 & 48 & 191 & 95 & 130 & 103 & \textbf{45} & 74 \\
yfitu & 47 & 45 & 48 & 69 & 95 & 130 & 103 & \textbf{43} & 74 \\
zecevic3 & 10 & 23 & \textbf{9} & 22 & 11 & 38 & 19 & 12 & 19 \\
zecevic4 & 7 & 13 & \textbf{6} & 10 & 8 & 22 & 12 & 15 & 16 \\
zigzag & 12 & 28 & \textbf{11} & 23 & 23 & 218 & 43 & 29 & 74 \\
zy2 & 6 & 14 & 5 & 10 & 9 & 43 & 9 & 14 & \textbf{4} \\
\hline
\end{longtable}
\end{scriptsize}

%% file: sections/solvers.tex

\begin{landscape}
\begin{table}[htbp]
\centering
\setlength{\tabcolsep}{5pt}
\renewcommand{\arraystretch}{1.5}
\caption{Description of state-of-the-art solvers within the proposed \revision{unifying} framework.}
\begin{tabular}{|p{1.05in}||p{1.1in}|p{1in}|p{0.78in}|p{0.8in}|p{0.6in}|p{0.74in}||p{0.4in}|}
\hline
Solver   & \cellcolor{constraint_relaxation_color} Constraint relaxation strategy & \cellcolor{inequality_method_color} Inequality handling method & \cellcolor{strategy_color} Globalization strategy   & \cellcolor{mechanism_color} Globalization mechanism   & \cellcolor{subproblem_color} Hessian model & \cellcolor{subproblem_color} Inertia correction strategy & Method \\
\hline
ALGLIB MinNLC     & $\ell_1$ relaxation         & IC & filter method     & trust region   & \revision{QN} & \revision{none} & SQP \\ \hline
CONOPT            & feasible step method        & IC & objective merit    & line search   & QN & \revision{none} & GRG \\ \hline
FICO XSLP         & $\ell_1$ relaxation         & IC & --               & $\ell_\infty$ trust region   & none & none & SLP \\ \hline
filterSQP         & feasibility restoration     & IC & filter method         & $\ell_\infty$ trust region   & exact & none & SQP \\ \hline
IPOPT             & feasibility restoration     & primal-dual IPM & filter method      & line search   & exact / QN & primal-dual & IPM \\ \hline
KNITRO-ASM        & $\ell_1$ relaxation         & EC & $\ell_1$ merit   & trust region   & exact & ? & SLQP \\ \hline
KNITRO-IPM        & \revision{step decomposition} & primal-dual IPM & $\ell_2$ merit   & line search + trust region   & exact & \revision{none} & IPM \\ \hline
LANCELOT          & bound-constrained AL        & EC & AL merit   & $\ell_2$ / $\ell_\infty$ trust region   & exact & none & AL \\ \hline
LOQO              & $\ell_2^2$ relaxation       & primal-dual IPM & $\ell_2^2$ merit   & line search   & exact & primal & IPM \\ \hline
MINOS             & linearly constrained AL + $\ell_1$ relaxation   & IC & forcing \newline sequences   & line search   & QN & \revision{none} & GRG \\ \hline
NAG e04uc/e04wdc  & $\ell_1$ relaxation         & IC & AL merit   & line search   & QN & \revision{none} & SQP \\ \hline
NLPQL             & right-hand-side \newline relaxation  & IC & AL merit   & line search & QN & \revision{none} & SQP \\ \hline
SLSQP             & right-hand-side \newline relaxation  & IC & $\ell_1$ merit   & line search   & QN & \revision{none} & SQP \\ \hline
SNOPT             & $\ell_1$ relaxation         & IC & AL merit   & line search   & QN & \revision{none} & SQP \\ \hline
SQuID             & $\ell_1$ relaxation         & IC & $\ell_1$ merit   & line search   & exact & primal & S$\ell_1$QP \\ \hline
WORHP             & right-hand-side \newline relaxation  & IC & $\ell_1$ merit / filter method   & line search & QN & \revision{none} & SQP \\ \hline
\end{tabular}
\label{tab:state-of-the-art-solvers}
\end{table}
\end{landscape}

%% file: sections/uml.tex
\begin{landscape}
\begin{figure}[h!]
\centering
\sffamily
\scalebox{0.8}{\input{figures/UML/uno_uml_paper}}
\caption{Uno's UML diagram.}
\label{fig:uno-uml}
\end{figure}
\end{landscape}

%% file: figures/UML/uno_uml_paper.tex
\begin{tikzpicture}
\tikzset{mechanism/.style={fill=mechanism_color}}
\tikzset{constraint_relaxation/.style={fill=constraint_relaxation_color}}
\tikzset{inequality_handling/.style={fill=inequality_method_color}}
\tikzset{strategy/.style={fill=strategy_color}}
\tikzset{subproblem/.style={fill=subproblem_color}}
\tikzset{subclass/.style={opacity=1}}

\umlclass[x=0, y=2.2]{Uno}{
   globalization\_mechanism
}{
	solve()
}

\umlabstract[mechanism]{GlobalizationMechanism}{
	constraint\_relaxation\_strategy
}{ 
	compute\_acceptable\_iterate()
}
\umlsimpleclass[mechanism, subclass, x=-5.4]{BacktrackingLineSearch}
\umlsimpleclass[mechanism, subclass, x=5.35]{TrustRegionStrategy}

\umlabstract[constraint_relaxation, x=0, y=-2.4]{ConstraintRelaxationStrategy}{
	inequality\_handling\_method
}{ 
	compute\_feasible\_direction() \\
	is\_acceptable()
}
\umlsimpleclass[constraint_relaxation, subclass, x=-5.35, y=-2.4]{FeasibilityRestoration}
\umlsimpleclass[constraint_relaxation, subclass, x=5.25, y=-2.4]{l1Relaxation}

\umlabstract[inequality_handling, x=-4.5, y=-4.7]{InequalityHandlingMethod}{
   subproblem
}{
   compute\_feasible\_direction()
}
\umlsimpleclass[inequality_handling, subclass, x=-9.5, y=-4.5]{InequalityConstrainedMethod}
\umlsimpleclass[inequality_handling, subclass, x=-9.5, y=-5.3]{InteriorPointMethod}

\umlabstract[subproblem, x=-6.5, y=-7.1]{Subproblem}{
    hessian\_model \\
    inertia\_correction\_strategy
}{
}

\umlabstract[x=-2.2, y=-6.9]{SubproblemSolver}{
}{ 
   solve()
}

\umlabstract[subproblem, x=-8, y=-9]{HessianModel}{
}{
}
\umlsimpleclass[subproblem, subclass, x=-11.5, y=-8.5]{ExactHessian}
\umlsimpleclass[subproblem, subclass, x=-11.5, y=-9.3]{IdentityHessian}
\umlsimpleclass[subproblem, subclass, x=-11.5, y=-10.1]{ZeroHessian}

\umlabstract[subproblem, x=-4, y=-9]{InertiaCorrectionStrategy}{
}{
}
\umlsimpleclass[subproblem, subclass, x=1, y=-8.5]{PrimalInertiaCorrection}
\umlsimpleclass[subproblem, subclass, x=1, y=-9.3]{PrimalDualInertiaCorrection}
\umlsimpleclass[subproblem, subclass, x=1, y=-10.1]{NoInertiaCorrection}

\umlabstract[strategy, x=4, y=-4.7]{GlobalizationStrategy}{
}{ 
   check\_acceptance()
}
\umlsimpleclass[strategy, subclass, x=2, y=-6.3]{l1MeritFunction}
\umlsimpleclass[strategy, subclass, x=6, y=-6.3]{FilterMethod}

\umlcompo{GlobalizationMechanism}{Uno}
\umlimpl[subclass]{BacktrackingLineSearch}{GlobalizationMechanism}
\umlimpl[subclass]{TrustRegionStrategy}{GlobalizationMechanism}

\umlcompo{ConstraintRelaxationStrategy}{GlobalizationMechanism}
\umlimpl[subclass]{FeasibilityRestoration}{ConstraintRelaxationStrategy}
\umlimpl[subclass]{l1Relaxation}{ConstraintRelaxationStrategy}

\umlcompo{InequalityHandlingMethod}{ConstraintRelaxationStrategy}
\umlimpl[subclass]{InequalityConstrainedMethod}{InequalityHandlingMethod}
\umlimpl[subclass]{InteriorPointMethod}{InequalityHandlingMethod}

\umlcompo{SubproblemSolver}{InequalityHandlingMethod}

\umlcompo{GlobalizationStrategy}{ConstraintRelaxationStrategy}
\umlimpl[subclass]{l1MeritFunction}{GlobalizationStrategy}
\umlimpl[subclass]{FilterMethod}{GlobalizationStrategy}

\umlcompo{Subproblem}{InequalityHandlingMethod}

\umlcompo{HessianModel}{Subproblem}
\umlimpl[subclass]{ExactHessian}{HessianModel}
\umlimpl[subclass]{IdentityHessian}{HessianModel}
\umlimpl[subclass]{ZeroHessian}{HessianModel}

\umlcompo{InertiaCorrectionStrategy}{Subproblem}
\umlimpl[subclass]{PrimalInertiaCorrection}{InertiaCorrectionStrategy}
\umlimpl[subclass]{PrimalDualInertiaCorrection}{InertiaCorrectionStrategy}
\umlimpl[subclass]{NoInertiaCorrection}{InertiaCorrectionStrategy}
\end{tikzpicture}

%% file: sections/combining_filtersqp.tex
\begin{algorithm}[htbp!]
\caption{Uno: trust-region filter restoration SQP.}
\label{alg:uno-tr-filter-sqp}
\scriptsize
\SetAlgoVlined
\KwIn{initial primal-dual iterate $(x^{(0)}, \constraintmultipliers^{(0)}, \boundmultipliers^{(0)})$, initial trust-region radius $\Delta$}
$k \gets 0$ \;
\begin{constraint_relaxation}
$(\infeasibilitymeasure(x), \objectivemeasure{\pi}(x), \highlight[subproblem_color]{\auxiliarymeasure(x)}) \equaldef (\| c(x) \|_1, f(x), \highlight[subproblem_color]{0})$ \;
$(\Delta \infeasibilitymodel^{(k)}(\direction{x}), \Delta \objectivemodel{\pi}^{(k)}(\direction{x}), \highlight[subproblem_color]{\Delta \auxiliarymodel^{(k)}(\direction{x})}) \equaldef \Big( \| c^{(k)} \|_1 - \| c^{(k)} + (\nabla c^{(k)})^T \direction{x} \|_1, -(\nabla f^{(k)})^T \direction{x}, \allowbreak \highlight[subproblem_color]{0} \Big)$ \;
$\mathit{phase} \gets$ Optimality \;
\end{constraint_relaxation}
\begin{globalization_strategy}
$\filterobjectivemeasure \equaldef \highlight[constraint_relaxation_color]{\objectivemeasure{1}} + \highlight[subproblem_color]{\auxiliarymeasure}$ \;
$\Delta \filterobjectivemodel^{(k)} \equaldef \highlight[constraint_relaxation_color]{\Delta \objectivemodel{1}^{(k)}} + \highlight[subproblem_color]{\Delta \auxiliarymodel^{(k)}}$ \;
Initialize $\mathcal{F}$ \;
\end{globalization_strategy}

\Repeat{termination criteria are satisfied} {
	\begin{globalization_mechanism}
	Set inner iteration counter $l \gets 0$ \;
	Reset trust-region radius $\Delta^{(l)} \in [\underline{\Delta}, \overline{\Delta}]$ \;
	\Repeat{$(\trial{x}^{(k+1,l)}, \trial{\constraintmultipliers}^{(k+1,l)}, \trial{\boundmultipliers}^{(k+1,l)})$ is acceptable} {
		\begin{constraint_relaxation}
      \If{$\mathit{phase}$ = Optimality} {
         $(\direction{x}^{(k,l)}, \direction{\constraintmultipliers}^{(k,l)}, \direction{\boundmultipliers}^{(k,l)}) \gets$ solve $\highlight[subproblem_color]{QP^{(k)}(\highlight[mechanism_color]{\Delta^{(l)}})}$ \;
         
         \If{$\highlight[subproblem_color]{QP^{(k)}(\highlight[mechanism_color]{\Delta^{(l)}})}$ infeasible} {
             $\mathit{phase} \gets$ Restoration \;
             $\constraintmultipliers^{(k)} \gets 0$ \;
             \begin{globalization_strategy}
             $\mathcal{F} \gets \mathcal{F} \cup \left\{ \left( \infeasibilitymeasure(x^{(k)}), \filterobjectivemeasure(x^{(k)}) \right) \right\}$ \;
             \end{globalization_strategy}
         }
      }
      \If{$\mathit{phase}$ = Restoration} {
			$(\direction{x}^{(k,l)}, \direction{\constraintmultipliers}^{(k,l)}, \direction{\boundmultipliers}^{(k,l)}) \gets$ solve $\highlight[subproblem_color]{FQP^{(k)}(\highlight[mechanism_color]{\Delta^{(l)}})}$ starting from $\direction{x}^{(k,l)}$ \;
		}
		\end{constraint_relaxation}
		
		Assemble trial iterate $(\trial{x}^{(k+1,l)}, \trial{\constraintmultipliers}^{(k+1,l)}, \trial{\boundmultipliers}^{(k+1,l)}) \equaldef (x^{(k)}, \constraintmultipliers^{(k)}, \boundmultipliers^{(k)}) + (\direction{x}^{(k,l)}, \direction{\constraintmultipliers}^{(k,l)}, \direction{\boundmultipliers}^{(k,l)})$ \;
      Reset the bound multipliers corresponding to the active trust region \;
      $\mathit{acceptable} \gets \mathit{false}$ \;
      
      \begin{constraint_relaxation}
		\eIf{$\|\direction{x}^{(k,l)}\| = 0$} {
			$\mathit{acceptable} \gets \mathit{true}$ \;
		}{
			\If{$\mathit{phase}$ = Restoration and $\highlight[subproblem_color]{QP^{(k)}(\highlight[mechanism_color]{\Delta^{(l)}})}$ feasible and $\highlight[strategy_color]{\infeasibilitymeasure(\trial{x}^{(k+1, l)}) ~< \eta_\mathit{min}(\mathcal{F})}$} {
				$\mathit{phase} \gets$ Optimality \;
            $\highlight[strategy_color]{\mathcal{F} \gets \mathcal{F} \cup \left\{ \left( \infeasibilitymeasure(x^{(k)}), \filterobjectivemeasure(x^{(k)}) \right) \right\}}$ \;
			}

         \begin{globalization_strategy}
         \uIf{$\mathit{phase}$ = Restoration} {
            \If{$\Delta \infeasibilitymeasure(x^{(k)}) - \infeasibilitymeasure(\trial{x}^{(k+1,l)}) \geq \sigma \Delta \infeasibilitymodel^{(k)}(\direction{x}^{(k,l)})$} {
               $\mathit{acceptable} \gets \mathit{true}$ \;
            }
         }
			\ElseIf{$\trial{x}^{(k+1,l)}$ acceptable to $\mathcal{F}$ and improves upon $x^{(k)}$} {
				\eIf{$\Delta \filterobjectivemodel^{(k)}(\direction{x}^{(k,l)}) \ge \delta \infeasibilitymeasure(x^{(k)})^2$} {
               \If{$\filterobjectivemeasure(x^{(k)}) - \filterobjectivemeasure(\trial{x}^{(k+1,l)}) \geq \sigma \Delta \filterobjectivemodel^{(k)}(\direction{x}^{(k,l)})$} {
                  $\mathit{acceptable} \gets \mathit{true}$ (f-type) \;
               }
				} {
					$\mathit{acceptable} \gets \mathit{true}$ (h-type) \;
               $\mathcal{F} \gets \mathcal{F} \cup \left\{ \left( \infeasibilitymeasure(x^{(k)}), \filterobjectivemeasure(x^{(k)}) \right) \right\}$ \;
            }
			}
         \end{globalization_strategy}
		}
      \end{constraint_relaxation}
		
      \eIf{$acceptable$} {
			\If{trust region is active at $\direction{x}^{(k,l)}$} {
				Increase radius $\Delta^{(l)}$ \;
			}
		}{
			Decrease radius $\Delta^{(l)}$ \;
		}
		$l \gets l+1$ \;
	}
	Update $(x^{(k+1)}, \constraintmultipliers^{(k+1)}, \boundmultipliers^{(k+1)}) \gets (\trial{x}^{(k+1,l)}, \trial{\constraintmultipliers}^{(k+1,l)}, \trial{\boundmultipliers}^{(k+1,l)})$ \;
	\end{globalization_mechanism}
	$k \gets k+1$ \;
}
\Return $(x^{(k)}, \constraintmultipliers^{(k)}, \boundmultipliers^{(k)})$
\end{algorithm}

%% file: sections/combining_ipopt.tex
\begin{algorithm}[htbp!]
\caption{Uno: line-search filter restoration interior-point method.}
\label{alg:uno-ls-filter-interior-point}
\scriptsize
\SetAlgoVlined
\KwIn{initial primal-dual iterate $(x^{(0)}, \constraintmultipliers^{(0)}, \boundmultipliers^{(0)})$, initial barrier parameter $\mu > 0$}
$k \gets 0$ \;
\begin{constraint_relaxation}
$(\infeasibilitymeasure(x), \objectivemeasure{\pi}(x), \highlight[subproblem_color]{\auxiliarymeasure(x)}) \equaldef (\| c(x) \|_1, f(x), \highlight[subproblem_color]{-\mu \log(X) e})$ \;
$(\Delta \infeasibilitymodel^{(k)}(\direction{x}), \Delta \objectivemodel{\pi}^{(k)}(\direction{x}), \highlight[subproblem_color]{\Delta \auxiliarymodel^{(k)}(\direction{x})}) \equaldef \Big( \| c^{(k)} \|_1 - \| c^{(k)} + (\nabla c^{(k)})^T \direction{x} \|_1, -(\nabla f^{(k)})^T \direction{x}, \allowbreak \highlight[subproblem_color]{\mu (X^{(k)})^{-1} e^T \direction{x}} \Big)$ \;
$\mathit{phase} \gets$ Optimality \;
\end{constraint_relaxation}
\begin{globalization_strategy}
$\filterobjectivemeasure \equaldef \highlight[constraint_relaxation_color]{\objectivemeasure{1}} + \highlight[subproblem_color]{\auxiliarymeasure}$ \;
$\Delta \filterobjectivemodel^{(k)} \equaldef \highlight[constraint_relaxation_color]{\Delta \objectivemodel{1}^{(k)}} + \highlight[subproblem_color]{\Delta \auxiliarymodel^{(k)}}$ \;
Initialize $\mathcal{F}$ \;
\end{globalization_strategy}

\Repeat{termination criteria are satisfied} {
	\begin{globalization_mechanism}
   \begin{constraint_relaxation}
   \begin{subproblem}
   Possibly update the barrier parameter $\mu$ \;
   \end{subproblem}
   $(\direction{x}^{(k)}, \direction{\constraintmultipliers}^{(k)}, \direction{\boundmultipliers}^{(k)}) \gets$ solve
   $\begin{cases}
   \highlight[subproblem_color]{IPSP_{\barrierparameter}^{(k)}}   & \text{ if } \mathit{phase} \text{ = Optimality} \\
   \highlight[subproblem_color]{FIPSP_{\barrierparameter}^{(k)}}  & \text{ if } \mathit{phase} \text{ = Restoration}
   \end{cases}$ \;
   \begin{subproblem}
   Scale primal-dual direction according to  (Eq.~\ref{eq:fraction-to-boundary}) \;
   \end{subproblem}
   \end{constraint_relaxation}
   
   $\alpha^{(0)} \gets 1$ \;
   Set inner iteration counter $l \gets 0$ \;
   
	\Repeat{$(\trial{x}^{(k+1, l)}, \trial{\constraintmultipliers}^{(k+1, l)}, \trial{\boundmultipliers}^{(k+1, l)})$ is acceptable} {
		Assemble trial iterate $(\trial{x}^{(k+1,l)}, \trial{\constraintmultipliers}^{(k+1,l)}, \trial{\boundmultipliers}^{(k+1,l)}) \equaldef (x^{(k)}, \constraintmultipliers^{(k)}, \boundmultipliers^{(k)}) + (\alpha^{(l)} \direction{x}^{(k)}, \alpha^{(l)} \direction{\constraintmultipliers}^{(k)}, \direction{\boundmultipliers}^{(k)})$ \;
      $\mathit{acceptable} \gets \mathit{false}$ \;

      \begin{constraint_relaxation}
		\eIf{$\|\direction{x}^{(k)}\| = 0$} {
			$\mathit{acceptable} \gets \mathit{true}$ \;
		}{
         \begin{globalization_strategy}
			\uIf{$\mathit{phase}$ = Restoration} {
            \If{$\Delta \infeasibilitymeasure(x^{(k)}) - \infeasibilitymeasure(\trial{x}^{(k+1,l)}) \geq \sigma \Delta \infeasibilitymodel^{(k)}(\direction{x}^{(k,l)})$} {
               $\mathit{acceptable} \gets \mathit{true}$ \;
            }
         }
			\ElseIf{$\trial{x}^{(k+1,l)}$ acceptable to $\mathcal{F}$} {
				\uIf{$\infeasibilitymeasure(x^{(k)}) \le \theta_\mathit{min}$ and $0 < \Delta \filterobjectivemodel^{(k)}(\highlight[mechanism_color]{\alpha^{(l)}} \direction{x}^{(k)})$ and $\Delta \filterobjectivemodel^{(k)}(\highlight[mechanism_color]{\alpha^{(l)}} \direction{x}^{(k)}) \ge \delta \infeasibilitymeasure(x^{(k)})^2$} {
               \eIf{$\filterobjectivemeasure(x^{(k)}) - \filterobjectivemeasure(\trial{x}^{(k+1,l)}) \geq \sigma \Delta \filterobjectivemodel^{(k)}(\highlight[mechanism_color]{\alpha^{(l)}} \direction{x}^{(k)})$} {
                  $\mathit{acceptable} \gets \mathit{true}$ \;
               } {
                  $\mathcal{F} \gets \mathcal{F} \cup \left\{ \left( \infeasibilitymeasure(x^{(k)}), \filterobjectivemeasure(x^{(k)}) \right) \right\}$ \;
               }
				}
            \ElseIf{$\trial{x}^{(k+1,l)}$ improves upon $x^{(k)}$} {
               $\mathit{acceptable} \gets \mathit{true}$ \;
            }
            \If{$\neg \left( 0 < \Delta \filterobjectivemodel^{(k)}(\highlight[mechanism_color]{\alpha^{(l)}} \direction{x}^{(k)}) \text{ and } \Delta \filterobjectivemodel^{(k)}(\highlight[mechanism_color]{\alpha^{(l)}} \direction{x}^{(k)}) \ge \delta \infeasibilitymeasure(x^{(k)})^2 \right)$} {
               $\mathcal{F} \gets \mathcal{F} \cup \left\{ \left( \infeasibilitymeasure(x^{(k)}), \filterobjectivemeasure(x^{(k)}) \right) \right\}$ \;
            }
			}
         \end{globalization_strategy}
         \If{$\mathit{acceptable}$ and $\mathit{phase}$ = Restoration and $\highlight[strategy_color]{\trial{x}^{(k+1, l)} \text{ acceptable to } \mathcal{F} \text{ and } \infeasibilitymeasure(\trial{x}^{(k+1, l)}) \le \kappa \infeasibilitymeasure(x^{(k)})}$} {
            $\mathit{phase} \gets$ Optimality \;
            $\highlight[strategy_color]{\mathcal{F} \gets \mathcal{F} \cup \left\{ \left( \infeasibilitymeasure(x^{(k)}), \filterobjectivemeasure(x^{(k)}) \right) \right\}}$ \;
         }
		}
      \end{constraint_relaxation}
		
		\If{not $acceptable$} {
			Decrease step length $\alpha^{(l)}$ \;
		}
		$l \gets l+1$ \;
      \If{$\alpha^{(l)}$ too small} {
         \begin{constraint_relaxation}
         $\mathit{phase} \gets$ Restoration \;
          \begin{globalization_strategy}
          $\mathcal{F} \gets \mathcal{F} \cup \left\{ \left( \infeasibilitymeasure(x^{(k)}), \filterobjectivemeasure(x^{(k)}) \right) \right\}$ \;
          \end{globalization_strategy}
         $(\direction{x}^{(k)}, \direction{\constraintmultipliers}^{(k)}, \direction{\boundmultipliers}^{(k)}) \gets$ solve $\highlight[subproblem_color]{FIPSP_{\barrierparameter}^{(k)}}$ \;
         \begin{subproblem}
         Scale primal-dual direction according to  (Eq.~\ref{eq:fraction-to-boundary}) \;
         \end{subproblem}
         \end{constraint_relaxation}
      }
	}
	Update $(x^{(k+1)}, \constraintmultipliers^{(k+1)}, \boundmultipliers^{(k+1)}) \gets (\trial{x}^{(k+1)}, \trial{\constraintmultipliers}^{(k+1)}, \trial{\boundmultipliers}^{(k+1)})$ \;
	\end{globalization_mechanism}
	$k \gets k+1$ \;
}
\Return $(x^{(k)}, \constraintmultipliers^{(k)}, \boundmultipliers^{(k)})$
\end{algorithm}